\documentclass[12pt,reqno]{amsart}
\usepackage{amscd}
\usepackage{hyperref}
\usepackage{latexsym,amssymb}
\usepackage{enumerate}
\usepackage{color}
\usepackage{amsthm}
\usepackage{amsmath}
\usepackage{graphicx}
\newtheorem{corollary}{Corollary}[section]

\newtheorem{lemma}[corollary]{Lemma}
\newtheorem{proposition}[corollary]{Proposition}
\newtheorem{remark}[corollary]{Remark}

\newtheorem{theorem}[corollary]{Theorem}

\newcommand{\mylabel}[1]{\label{#1}
            \ifx\undefined\stillediting
            \else \fbox{$#1$}\fi }
\newcommand{\BE}{\begin{equation}}

\newcommand{\EEQ}{\end{equation}}
\newcommand{\rfb}[1]{\mbox{\rm
   (\ref{#1})}\ifx\undefined\stillediting\else:\fbox{$#1$}\fi}

\newfont{\Blackboard}{msbm10 scaled 1200}

\newfont{\roma}{cmr10 scaled 1200}

\newcommand{\mm}    {{\hbox{\hskip 0.5pt}}}

\newcommand{\bluff} {{\hbox{\raise 15pt \hbox{\mm}}}}
\newcommand{\e}      {{\varepsilon}}

 \allowdisplaybreaks \textwidth=16cm \textheight=23cm
\makeatletter
\def\section{\@startsection {section}{1}{\z@}{-3.5ex plus -1ex minus
    -.2ex}{2.3ex plus .2ex}{\large\bf}}
\makeatother
\def\be{\begin{equation}}
\def\ee{\end{equation}}

\begin{document}
\thispagestyle{empty}
\title{Null controllability of a cascade model in population dynamics}

\author{Bedr'Eddine Ainseba, Younes Echarroudi and Lahcen Maniar }
\thanks{ Institut de Math\'ematiques de Bordeaux, UMR-CNRS 5251, Universit\'{e} Bordeaux Segalen, 3 Place de la Victoire, 33076 Bordeaux Cedex, France,
e-mail: bedreddine.ainseba@u-bordeaux2.fr}
\thanks{Private university of Marrakesh, Km 13 Route d'Amizmiz, Marrakesh, Morocco, \\e-mail: yecharroudi@gmail.com,}
\thanks{D\'epartement de
Math\'ematiques,
Facult\'e des Sciences Semlalia, Laboratoire LMDP, UMMISCO (IRD-UPMC), B. P. 2390 Marrakech 40000, Maroc,
e-mail: maniar@uca.ma}


\begin{abstract}
In this paper, we are concerned
with the null controllability of a linear population dynamics cascade systems (or the so-called  prey-predator models) with two different dispersion
coefficients which degenerate in the boundary and with one control force. We develop first a  Carleman type
inequality for its adjoint system,  and then an observability inequality which allows us to deduce the
existence of a control acting on a subset of  the space domain  which steers both populations of a certain
age  to extinction in  a finite time.
\end{abstract}

\maketitle

\section{Introduction}\label{first-section-5}
We consider the coupled population cascade system
\begin{eqnarray}
  {{\partial y} \over {\partial t}} + {{\partial y} \over {\partial a}}-(k_{1}(x)y_{x})_{x}+\mu_{1}(t,a, x)y = \vartheta\chi_{\omega}   && \text{ in } Q, \label{470}\\
\nonumber {{\partial p} \over {\partial t}} + {{\partial p} \over {\partial a}}-(k_{2}(x)p_{x})_{x}+\mu_{2}(t,a, x)p+ \mu_{3}(t,a, x)y =0 && \text{ in } Q, \\
\nonumber  y(t,a, 1)=y(t,a, 0)=p(t,a, 1)=p(t,a, 0)=0  && \text{ on }(0,A)\times(0,T),\\
\nonumber y(0, a, x)=y_{0}(a, x);  p(0, a, x)=p_{0}(a, x) &&\text{ in }Q_{A},\\
\nonumber y(t, 0, x)=\int_{0}^{A} \beta_{1} (t,a, x)y (t,a, x)da  &&  \text{ in } Q_{T},\\
\nonumber p(t, 0, x)=\int_{0}^{A} \beta_{2} (t,a, x)p (t,a, x)da  && \text{ in } Q_{T},
\end{eqnarray}
where $Q=(0,T)\times(0,A)\times(0,1)$, $Q_{A}=(0,A)\times(0,1)$, $Q_{T}=(0,T)\times(0,1)$ and we will denote
 $q=(0,T)\times(0,A)\times \omega$. The system \eqref{470} models the dispersion of a gene in two given populations which are in interaction. In this case, $x$ represents the gene type and $y(t,a,x)$ and $p(t,a,x)$ as the distributions of individuals of age $a$ at time $t$ and of gene type $x$ of both populations. The parameters $\beta_{1}(t,a,x)$ (respectively
 $\beta_{2}(t,a,x)$),  $\mu_{1}(t,a,x)$ (respectively $\mu_{2}(t,a,x)$) are respectively the natural fertility
 and mortality rates of individuals of age $a$ at time $t$ and of gene type $x$ of the population whose
 distribution is $y$ (respectively $p$), $\mu_{3}$ can be interpreted as the interaction coefficient between two populations (cancer cells and healthy cells for instance)  which depends on $x$, $t$ and $a$,
 the subset $\omega$ is the region where a control $\vartheta$ is acting. Such a control corresponds to an external supply or to removal of individuals on the subdomain  $\omega$. Finally, $\int_{0}^{A} \beta_{1} (t,a, x)y (t,a, x)da$ and $\int_{0}^{A} \beta_{2} (t,a, x)p(t,a, x)da$ are
 the distributions of the newborns of the two populations that are of gene type $x$ at time $t$.\\\\ The control problems of \eqref{470} or in general of coupled systems take an intense interest and are widely investigated in many papers, among
them we find \cite{Ain5}, \cite{Haj1}, \cite{Zh2} and the references therein. In fact, in \cite{Ain5} the authors studied
a coupled reaction-diffusion equations describing interaction between a prey population and predator population.
The goal of this work was to look for a suitable control supported on a small spatial subdomain which guarantees
the stabilization of the predator population to zero. In \cite{Zh2}, the objective was different. More precisely,
the authors considered an age-dependent prey-predator system and they proved the existence and uniqueness for an
optimal control (called also "optimal effort") which gives the maximal harvest via the study of the optimal
harvesting problem associated to their coupled model.
\\However, the previous results were found in the case when
the diffusion coefficients are constants. This leads Ait Ben Hassi et al. in \cite{Haj1} to generalize the model
of \cite{Ain5} and investigate a semilinear parabolic cascade systems with two different diffusion coefficients
allowed to depend on the space variable and degenerate at the left boundary of the space domain. Moreover,
the purpose of this paper was to show the null controllability via a Carleman type inequality of the adjoint problem
of the associated linearized system using the results of \cite{Bouss} (or \cite{Can1}) and with the help of the Schauder fixed point theorem. On the other hand, a massive interest was given to the question of null controllability of the population dynamics models in the case of one equation both in the case without diffusion (see for example \cite{Marcheva}) and with diffusion (see for instance \cite{Ain4, Ain3, Ain2, Ain1, Traore} in the case of a constant diffusion coefficient). Recently, a more general case was investigated by B. Ainseba and al. in \cite{ech} and \cite{man}. Indeed, in \cite{ech} the authors allowed the dispersion coefficient to depend on the variable $x$ and verifies $k(0)=0$ (i.e, the coefficient of dispersion $k$ degenerates at 0) and they tried to obtain the null controllability  in such a situation with $\beta \in L^{\infty}$ basing on the work done in \cite{Bouss} for the degenerate heat equation to establish a new Carleman estimate for a suitable full adjoint system and afterwards his observability inequality. However, the main controllability result of \cite{ech} was shown under the condition $T\geq A$ (as in \cite{Marcheva}) and this constitutes a restrictiveness on the "optimality" of the control time $T$ since it means, for example, that for a pest population whose the maximal age $A$ may equal to a many days (may be many months or years) we need much time to bring the population to the zero equilibrium. In the same trend and to overcome the condition $T\geq A$, L. Maniar et al in \cite{man} suggested the fixed point technique implemented in \cite{Traore} and which requires that the fertility rate must belong to $C^{2}(Q)$ and consists briefly to demonstrate in a first time the null controllability for an intermediate system with a fertility function $b\in L^{2}(Q_{T})$ instead of $ \int_{0}^{A} \beta (t,a, x)y (t,a, x)da $ and to achieve the task via a Leray-Schauder theorem.\\ But up now, little is known about the null controllability question of population dynamics cascade systems both in degenerate and nondegenerate cases to our knowledge and the work done in this paper will address to such a control problem and it will be a generalization of the results established in \cite{ech} and \cite{man}. More precisely, following the strategy of \cite{Haj1} we expect in this contribution to prove the null controllability of system \eqref{470} when $T\in (0, \delta)$ where $\delta\in(0,A)$ small enough in the case of one control force. That is, we show that for all $y_{0}, p_{0} \in L^{2}(Q_{A})$ and $\delta\in (0,A)$ small enough, there exists a control $\vartheta \in L^{2}(q)$ such that the associated solution of \eqref{470}
verifies
\begin{equation}\label{471}
\begin{cases}
&y(T,a, x)=0, \quad  \text{ a.e. in } (\delta, A)\times(0,1),\\
&p(T,a, x)=0, \quad   \text{ a.e in } (\delta, A)\times(0,1).
\end{cases}
\end{equation}
Such a result is gotten under the conditions that all the natural rates possess an $L^{\infty}-$regularity (see \eqref{474} beneath) and the dispersion coefficients are different and depend on the gene type with a degeneracy in the left hand side of its domain, i.e $k_{i}(0)=0; i=1,2$ (e.g $k_{i}=x^{\alpha}$, $\alpha>0$). In this case, we say that \eqref{470} is a degenerate population dynamics cascade system. Genetically speaking, such a property is natural since it means that if each population is not of a gene type, it can not be transmitted to its offspring.
\\The remainder of this paper is organized as follows: in Section \ref{second-section-5}, we  give the well-posedness result of system
\eqref{470} and we bring out a Carleman  inequality for an intermediate trivial adjoint system which helps us to prove the main  Carleman estimate for the full adjoint model. With the aid of this inequality,
we establish in Section \ref{third-section-5}  the observability inequality and show the main  result of the null controllability of \eqref{470}. The last section takes the form of an appendix wherein we will give the proofs of some basic tools.
\section{Well-posedness and  Carleman estimates}\label{second-section-5}
\subsection{Well-posedness result}\label{second-section-5-1}
For this section and for the sequel, we assume that the dispersion coefficients  $k_{i}, i=1, 2$ satisfy the  hypotheses
\begin{equation}\label{473}
\left\{
\begin{array}{l}
k_{i} \in  C([0,1])\cap C^{1}((0,1]),\;\;
k_{i}>0 \text{ in }(0,1] \text{ and } k_{i}(0)=0,\\
\exists \gamma \in [0,1):  xk_{i}'(x)\leq\gamma k_{i}(x), \;  x \in [0,1].
\end{array}
\right.
\end{equation}
The last hypothesis on $k_{i}$ means in the case of $k(x)=x^{\alpha_{i}}$ that $0\leq \alpha_{i}<1$.
Similarly,  all results of this paper can be obtained  also in the case of $1\leq \alpha_{i}<2$ taking,
instead of Dirichlet condition on $x=0$, the Newmann condition $(k_{i}(x)u_x)(0)=0$.
On the other hand, we assume that  the  rates $\mu_{1}$, $\mu_{2}$, $\mu_{3}$, $\beta_{1}$ and $\beta_{2}$ verify
\begin{equation}\label{474}
\begin{cases}
 \mu_{1}, \mu_{2}, \mu_{3}, \beta_{1}, \beta_{2} \in L^{\infty}(Q),\quad
 \mu_{1}, \mu_{2}, \mu_{3}, \beta_{1}, \beta_{2}\geq0  \text{ a.e in }Q,\\
\beta_{i}(.,0,.)\equiv0 \text{ a.e. in } (0,T)\times(0,1), \quad \text{ for } i=1,2.\\
\end{cases}
\end{equation}
The third assumption in \eqref{474} on the fertility rates $\beta_{1}$ and $\beta_{2}$ is natural since the newborns are not fertile.
 \\As in \cite{man}, we discuss the well-posedness of  \eqref{470} by introducing the weighted
spaces $H^{1}_{k_{i}}(0, 1)$ and $H^{2}_{k_{i}}(0, 1)$ defined by
\begin{equation*}\label{578}
\left\{
\begin{array}{l}
H^{1}_{k_{i}}(0, 1):=\{u \in L^{2}(0,1):u \text{ is abs.  cont. in } [0,1]:  \sqrt{k_{i}}u_{x} \in L^{2}(0,1), u(1)=u(0)=0\},\\
H^2_{k_{i}} (0, 1):=\Big\{ u \in H^1_{k}(0, 1) \,: \, k_{i}(x)u_x
\in H^1(0, 1)\Big\},
\end{array}
\right.
\end{equation*}
endowed respectively with the norms
\begin{equation*}\label{579}
\left\{
\begin{array}{l}
\|u\|^{2}_{H^{1}_{k_{i}}(0, 1)}
:=\|u\|^{2}_{L^{2}(0,1)}+ \|\sqrt{k_{i}}u_{x}\|^{2}_{L^{2}(0,1)}, \quad u \in H^{1}_{k_{i}}(0, 1),\\
\|u\|^2_{H^2_{k_{i}}} := \|u\|^2 _{H^1_{k_{i}}(0, 1)}
+ \|(k_{i}(x)u_x)_x\|^2_{ L^2(0,1)}, \quad u \in H^{2}_{k_{i}}(0, 1),
\end{array}
\right.
\end{equation*}
with $i=1,2$ (see \cite{Haj1}, \cite{Bouss}, \cite{Can1} or the references therein for the properties of such a spaces).
 We recall from
\cite {cmp, Can1} that the  operators
$
C_{i}u := (k_{i}(x)u_x)_x,\, \ u \in
D(C_{i}) = H^2_{k_{i}}(0, 1), i=1,2
$
are closed
self-adjoint and  negative  with dense domains in $L^2(0, 1)$.
\\On the other hand, in the Hilbert space $\mathbb{H}=(L^{2}((0, A)\times(0, 1)))^{2}$, the system \eqref{470} can be rewritten abstractly as an inhomogeneous Cauchy problem in the following way
\begin{equation}
\begin{cases}\
X^{'}(t)=\mathbb{A}X(t)+B(t)X(t)+f(t),\\
X(0)=\left(
       \begin{array}{c}
         y_{0} \\
         p_{0}\\
       \end{array}
     \right),
\end{cases}
\end{equation}
where $X(t)=\left(
              \begin{array}{c}
                y(t) \\
                p(t)\\
              \end{array}
            \right)
$, $\mathbb{A}=\left(
                 \begin{array}{cc}
                   \mathcal{A}_{1}& 0\\
                   0 & \mathcal{A}_{2}\\
                 \end{array}
               \right)
$; $D(\mathbb{A})=D(\mathcal{A}_{1})\times D(\mathcal{A}_{2})$, \\$f(t)=\left(
                                                                \begin{array}{c}
                                                                  \vartheta(t,.,\cdot)\chi_{\omega}(.) \\
                                                                  0
                                                                \end{array}
                                                              \right)
$, $B(t)=\left(
           \begin{array}{cc}
             M_{\mu_{1}(t)} & 0 \\
              M_{\mu_{3}(t)} &  M_{\mu_{2}(t)} \\
           \end{array}
         \right)
$, where $ M_{\mu_{j}(t)}w=-\mu_{j}(t)w$, the operators $\mathcal{A}_{1}:L^{2}((0, A)\times(0, 1))\rightarrow L^{2}((0, A)\times(0, 1))$ and $\mathcal{A}_{2}:L^{2}((0, A)\times(0, 1))\rightarrow L^{2}((0, A)\times(0, 1))$ are defined respectively by:
\begin{equation}
\begin{cases}
\mathcal{A}_{1}\theta(a,x)=-{{\partial \theta} \over {\partial a}}+(k_{1}(x)\theta_{x})_{x}, \hspace{0,25cm} \forall \theta \in D(\mathcal{A}_{1}), \\
D(\mathcal{A}_{1})=\{\theta(a,x): \theta, \mathcal{A}_{1}\theta\in L^{2}((0, A)\times(0, 1)), \theta(a,0)=\theta(a,1)=0,
\theta(0,x)=\int_{0}^{A} \beta_{1} (a, x)\theta(a, x)da\},
\end{cases}
\end{equation}
and
\begin{equation}
\begin{cases}
\mathcal{A}_{2}\theta(a,x)=-{{\partial \theta} \over {\partial a}}+(k_{2}(x)\theta_{x})_{x}, \hspace{0,25cm} \forall \theta \in D(\mathcal{A}_{2}), \\
D(\mathcal{A}_{2})=\{\theta(a,x): \theta, \mathcal{A}_{2}\theta\in L^{2}((0, A)\times(0, 1)), \theta(a,0)=\theta(a,1)=0,
\theta(0,x)=\int_{0}^{A} \beta_{2} (a, x)\theta(a, x)da\}.
\end{cases}
\end{equation}
It is well-known, from \cite{web} and the references therein that the operators $\mathcal{A}_{1}$ and $\mathcal{A}_{2}$ defined above generate a $C_{0}-$semigroups. On the other hand, one can see that the operator $\mathbb{A}$ is diagonal and $B(t)$ is a bounded perturbation. Therefore, the following well-posedness result holds (see for instance \cite{Haj1} for a similar result of cascade parabolic equations).
\begin{theorem}\label{475}
$i)$ The operator $\mathbb{A}$ generates a $C_{0}-$semigroup.\\
$ii)$ Under the assumptions \eqref{473} and \eqref{474} and for all $\vartheta \in L^{2}(Q)$ and
$(y_{0}, p_{0}) \in (L^{2}(Q_{A}))^{2}$, the system \eqref{470} admits a unique solution $(y,p)$.
This solution belongs to $E:=C([0, T], (L^{2}((0, A)\times(0, 1)))^{2})\cap C([0, A], (L^{2}((0, T)\times(0, 1)))^{2}) \cap L^{2}((0, T)\times(0, A), H^{1}_{k_{1}}(0, 1)\times H^{1}_{k_{2}}(0, 1) )$.
Moreover, the solution of \eqref{470} satisfies the following \\inequality
\begin{eqnarray}\label{476}
 \nonumber & &
\sup_{t\in [0,T]}\|(y(t),p(t))\|^{2}_{L^{2}(Q_{A})\times L^{2}(Q_{A})}+ \sup_{a\in [0,A]}\|(y(a),p(a))\|^{2}_{L^{2}(Q_{T})\times L^{2}(Q_{T})}\\ \nonumber & &+\int_{0}^{1}\int_{0}^{A}\int_{0}^{T}
((\sqrt{k_{1}}y_x)^{2}+(\sqrt{k_{2}}p_x)^2)dtdadx\\  & &\leq C\left(\int_{q}\vartheta^{2}dtdadx+
\|(y_{0}, p_{0})\|^{2}_{L^{2}(Q_{A})\times L^{2}(Q_{A})}\right).
\end{eqnarray}
\end{theorem}
\subsection{Carleman inequality results}\label{second-section-5-2}

In this paragraph, we show a Carleman type inequality for the following adjoint system of \eqref{470} \begin{eqnarray}
{{\partial u} \over {\partial t}} + {{\partial u} \over {\partial a}}+(k_{1}(x)u_{x})_{x}-\mu_{1}(t,a, x)u- \mu_{3}(t,a, x)v =-\beta_{1}(t,a,x)u(t,0,x) && \text{ in } Q, \label{477}\\
{{\partial v} \over {\partial t}} + {{\partial v} \over {\partial a}}+(k_{2}(x)v_{x})_{x}-\mu_{2}(t,a, x)v =-\beta_{2}(t,a,x)v(t,0,x) &&\text{ in } Q,\nonumber\\
 u(t,a, 1)=u(t,a, 0)=v(t,a, 1)=v(t,a, 0)=0 && \text{ on } (0,T)\times(0,A), \nonumber \\
\nonumber  u(T,a, x)=u_{T}(a, x) && \text{ in } Q_{A}, \\
\nonumber  v(T,a, x)=v_{T}(a, x) && \text{ in } Q_{A},\\
\nonumber  u(t,A, x)=v(t,A, x)=0 && \text{ in } Q_{T}.
\end{eqnarray}
To do this, we prove firstly the Carleman estimate  for the following intermediate system
\begin{eqnarray}
{{\partial u} \over {\partial t}} + {{\partial u} \over {\partial a}}+(k_{1}(x)u_{x})_{x}-\mu_{1}(t,a, x)u- \mu_{3}(t,a, x)v =h_{1} && \text{ in } Q, \label{478}\\
{{\partial v} \over {\partial t}} + {{\partial v} \over {\partial a}}+(k_{2}(x)v_{x})_{x}-\mu_{2}(t,a, x)v=h_{2} && \text{ in } Q,\nonumber\\
 u(t,a, 1)=u(t,a, 0)=v(t,a, 1)=v(t,a, 0)=0 && \text{ on } (0,T)\times(0,A), \nonumber \\
\nonumber  u(T,a, x)=u_{T}(a, x) && \text{ in } Q_{A},\\
\nonumber  v(T,a, x)=v_{T}(a, x) && \text{ in } Q_{A},\\
\nonumber  u(t,A, x)=v(t,A, x)=0 && \text{ in } Q_{T},
\end{eqnarray}
with $(u_{T}, v_{T}) \in (L^{2}(Q_{A}))^{2}$ and  $h_{1}, h_{2} \in L^{2}(Q)$. Such a system can be rewritten in the following way
\begin{eqnarray}
{{\partial u} \over {\partial t}} + {{\partial u} \over {\partial a}}+(k_{1}(x)u_{x})_{x}-\mu_{1}(t,a, x)u=h_{1}+\mu_{3}(t,a, x)v  && \text{ in } Q, \label{inter-1}\\
u(t,a, 1)=u(t,a, 0)=0 && \text{ on } (0,T)\times(0,A), \nonumber \\
\nonumber  u(T,a, x)=u_{T}(a, x) && \text{ in } Q_{A},\\
\nonumber  u(t,A, x)=0 && \text{ in } Q_{T},
\end{eqnarray}
where $v$ is the solution of
\begin{eqnarray}
{{\partial v} \over {\partial t}} + {{\partial v} \over {\partial a}}+(k_{2}(x)u_{x})_{x}-\mu_{2}(t,a, x)v=h_{2}  && \text{ in } Q, \label{inter-2}\\
v(t,a, 1)=v(t,a, 0)=0 && \text{ on } (0,T)\times(0,A), \nonumber \\
\nonumber  v(T,a, x)=v_{T}(a, x) && \text{ in } Q_{A},\\
\nonumber  v(t,A, x)=0 && \text{ in } Q_{T}.
\end{eqnarray}
Classically, the proof of such a kind of estimates is based tightly on the choice of the so-called weight functions. In our case, these functions are set in the following way
\begin{equation}\label{479}
\left\{
\begin{array}{l}
\varphi_{i}(t,a, x):=\Theta  (t, a)\psi_{i}(x), i=1,2,
\\
\Theta  (t, a):= \displaystyle\frac{1}{(t(T-t))^{4}a^{4}},\\
\psi_{i}(x):= \lambda_{i}\left(\int_{0}^{x}\frac{r}{k_{i}(r)} dr-d_{i}\right),\\
\phi(t,a,x)=\Theta(a,t)e^{\kappa\sigma(x)},
\Phi(t,a,x)=\Theta(a,t)\Psi(x), \Psi(x)=e^{\kappa\sigma(x)}-e^{2\kappa\|\sigma\|_{\infty}},
\end{array}
\right.
\end{equation}
where $\sigma$ is the  function given by \begin{equation}\label{568}
\left\{
\begin{array}{l}
\sigma \in C^{2}([0,1]),
\sigma(x)>0 \text{ in } (0,1), \sigma(0)=\sigma(1)=0, \\
\sigma_{x}(x)\neq0 \text{ in } [0,1]\backslash \omega_{0},
\end{array}
\right.
\end{equation}
$ \omega_{0}\Subset\omega $ is an open subset. The existence of this function is proved in
\cite[Lemma 1.1]{Fursikov}. $\lambda_{i}$, $d_{i}$ for $i=1,2$ and $\kappa$ are supposed to verify following assumptions
\begin{equation}\label{parameters}
\left\{
\begin{array}{l}
d_{1}> \frac {1}{k_{1}(1)(2-\gamma)}, \frac{\lambda_{1}}{\lambda_{2}}\geq\frac{d_{2}}{d_{1}-\int_{0}^{1}\frac{r}{k_{1}(r)} dr},\\
\kappa\geq\frac{4\ln2}{\|\sigma\|_{\infty}}, d_{2}\geq \frac {5}{k_{2}(1)(2-\gamma)},
\end{array}
\right.
\end{equation}
with $\lambda_{2}\in I=[\frac{k_{2}(1)(2-\gamma)(e^{2\kappa\|\sigma\|_{\infty}}-1)}{d_{2}k_{2}(1)(2-\gamma)-1}, \frac{4(e^{2\kappa\|\sigma\|_{\infty}}-e^{\kappa\|\sigma\|_{\infty}})}{3d_{2}})$ which can be shown not empty (see Lemma \ref{Intervalle} in the appendix). On other hand,  in the light of the first and the fourth conditions in
\eqref{parameters} on $d_{1}$ and $d_{2}$, one can observe that $\psi_{i}(x)<0$ for  all
 $x \in [0,1]$, and $\Theta  (t, a)\rightarrow +\infty$ as $t\rightarrow 0^{+}, T^{-}$ and $a\rightarrow 0^{+}$.

Now, we state the first result of this section which is the intermediate Carleman estimate satisfied by solution of system \eqref{478}.
\begin{theorem}\label{480}
Assume that $k_{i}$ satisfy the hypotheses \eqref{473} and let $A>0$ and $T>0$ be given.
Then, there exist two positive constants $C$ and $s_{0}$, such that every solution $(u,v)$ of \eqref{478} satisfies,
for all $s\geq s_{0}$, the following inequality
\begin{eqnarray}\label{omega-estimate}
& &\int_{Q} \left(s^{3}\Theta^{3}\frac{x^{2}}{k_{1}(x)}u^{2}+s\Theta k_{1}(x)u_{x}^{2}\right)e^{2s\varphi_{1}} dtdadx
+ \int_{Q} \left(s^{3}\Theta^{3}\frac{x^{2}}{k_{2}(x)}v^{2}+s\Theta k_{2}(x)v_{x}^{2}\right)e^{2s\varphi_{2}} dtdadx \nonumber\\& &\leq C\left(\int_{Q}(h_{1}^{2}+h_{2}^{2})e^{2s\Phi} dtdadx
+\int_{q}s^{3}\Theta^{3}(u^{2}+v^{2})e^{2s\Phi}dtdadx\right).
\end{eqnarray}
\end{theorem}
The proof of Theorem \ref{480} needs two basic results. These results are concerned with Carleman type inequalities in both cases degenerate and nondegenerate. The first one is stated in the following proposition
\begin{proposition}\label{15}
Consider the following system with $h \in L^{2}(Q)$, $\mu \in L^{\infty}(Q)$ and $k$ verifies the hypotheses \eqref{473}
 \begin{eqnarray}\label{17}
\label{eq.mdl.1}
  {{\partial u} \over {\partial t}} + {{\partial u} \over {\partial a}}+(k(x)u_{x})_{x}-\mu(t, a, x)u &=& h,\\
\nonumber  u(t, a, 1)=u(t, a, 0)&=&0, \\
\nonumber  u(T, a, x)&=&u_{T}(a, x),\\
\nonumber  u(t, A, x)&=&0.
\end{eqnarray}
 Then, there exist two positive constants
$C$ and $s_{0}$, such that every solution of \eqref{17} satisfies, for all $s\geq s_{0}$, the following inequality
\begin{eqnarray}\label{16}
& &s^{3}\int_{Q}  \Theta^{3}\frac{x^{2}}{k(x)}u^{2}e^{2s\varphi} dtdadx + s\int_{Q}  \Theta k(x)u_{x}^{2}e^{2s\varphi} dtdadx \\ \nonumber& &\leq C\left(\int_{Q}\mid h\mid^{2}e^{2s\varphi} dtdadx+sk(1) \int_{0}^{A}\int_{0}^{T}\Theta u_{x}^{2}(a, t, 1)e^{2s\varphi(a, t, 1)} dtda\right),
\end{eqnarray}
where $\varphi$ and $\Theta$ are the weight functions defined by
\begin{equation}\label{13}
\left\{
\begin{array}{l}
\varphi(t, a, x):=\Theta (t, a)\psi(x) \text{ with :}\\
\Theta (t, a):= \frac{1}{(t(T-t))^{4}a^{4}},\\
\psi(x):=c_{1}(\int_{0}^{x}\frac{r}{k(r)} dr-c_{2}).
\end{array}
\right.
\end{equation}
with $c_{2}> \frac {1}{k(1)(2-\gamma)}$, $c_{1}>0$ and $\gamma$ is the parameter defined by \eqref{473}.
\end{proposition}
For the proof of this proposition, we refer the reader to \cite[Proposition 3.1]{man}. The second result is the following
\begin{proposition}\label{Carl-nondegenerate}
Let us consider the following system
\begin{eqnarray}\label{564}
{{\partial z} \over {\partial t}} + {{\partial z} \over {\partial a}}+(k(x)z_{x})_{x}-c(t, a, x)z = h && \text{ in } Q_{b},\\
\nonumber  z(t, a, b_{1})=z(t, a, b_{2})=0 && \text{ on } (0,T)\times(0,A),
\end{eqnarray}
where $Q_{b}:=(0,T)\times(0,A)\times(b_{1},b_{2})$, $(b_{1},b_{2})\subset[0,1]$, $h \in L^{2}(Q_{b})$, $k \in C^{1}([0,1])$ is a strictly positive function and $c\in L^{\infty}(Q_{b})$. Then, there exist two positive constants $C$ and $s_{0}$, such that for any $s\geq s_{0}$, $z$ verifies the following
estimate
\begin{eqnarray}\label{570}
\nonumber & &
\int_{Q_{b}}(s^{3}\phi^{3}z^{2}+s\phi z_{x}^{2})e^{2s\Phi}dtdadx \\ & &
\leq C \left(\int_{Q_{b}}h^{2}e^{2s\Phi}dtdadx+\int_{\omega}\int_{0}^{A}\int_{0}^{T}s^{3}\phi^{3}z^{2}e^{2s\Phi}dtdadx\right),
\end{eqnarray}
where
$\phi$, $\Theta$ and $\Phi$ are defined by \eqref{479} and $\sigma$ by \eqref{568}.
\end{proposition}
For the proof of Proposition \ref{Carl-nondegenerate}, a careful computations allow us to adapt the same procedure of \cite[Lemma 2.1]{Ain3} to show \eqref{570} in case where $k$ is a positive general nondegenerate coefficient, with our weight function $\Theta(t,a)=\frac{1}{t^{4}(T-t)^{4}a^{4}}$ and the source term $h$.\\\\ Besides the two Propositions \ref{15} and \ref{Carl-nondegenerate}, we must bring out another important result
\begin{lemma}\label{interval-2}
Under assumptions \eqref{parameters}, the functions $\varphi_{1}$, $\varphi_{2}$ and $\Phi$
defined by \eqref{479} satisfy the following inequalities
\begin{equation}\label{interval-3}
\left\{
\begin{array}{l}
\varphi_{1}\leq\varphi_{2},\\
\frac{4}{3}\Phi<\varphi_{2}\leq\Phi.
\end{array}
\right.
\end{equation}
\end{lemma}
\proof
By the definitions of $\varphi_{1}$, $\varphi_{2}$ and $\Phi$ and taking into account that $\Theta$ is positive, showing
the results of \eqref{interval-3} is equivalent to show
\begin{eqnarray}\label{interval-4}
\left\{
\begin{array}{l}
\psi_{1}\leq\psi_{2},\\
\frac{4}{3}\Psi<\psi_{2}\leq\Psi.
\end{array}
\right.
\end{eqnarray}
The first inequality in \eqref{interval-4} is assured by the second assumption in \eqref{parameters} while the second one is deduced from  $\lambda_{2}\in I=[\frac{k_{2}(1)(2-\gamma)(e^{2\kappa\|\sigma\|_{\infty}}-1)}{d_{2}k_{2}(1)(2-\gamma)-1}, \frac{4(e^{2\kappa\|\sigma\|_{\infty}}-e^{\kappa\|\sigma\|_{\infty}})}{3d_{2}})$ and
this achieves the proof.
\endproof
Now, we can address the proof of Theorem \ref{480}.
\proof
Let us introduce the smooth cut-off function $\xi:\mathbb{R}\rightarrow \mathbb{R}$ defined as follows
\begin{equation}\label{481}
\begin{cases}
0\leq\xi(x)\leq1,&x \in \mathbb{R},\\
\xi(x)=1, &   x \in [0,\frac{2x_1+x_2}{3}],\\
\xi(x)=0, &x \in [\frac{x_1+2x_2}{3},1].
\end{cases}
\end{equation}
Let $u$ and $v$ be respectively the  solutions of \eqref{inter-1} and \eqref{inter-2}. Set $w:=\xi u$, $z:=\xi v$ and
put $\omega^{'}=(\frac{2x_1+x_2}{3},\frac{x_1+2x_2}{3})$.
Then, $(w,z)$ satisfies the following system
\begin{eqnarray}
{{\partial w} \over {\partial t}} + {{\partial w} \over {\partial a}}+(k_{1}(x)w_{x})_{x}-\mu_{1}(t,a, x)w= \mu_{3}(t,a, x)z+\xi h_{1}+(k_{1}\xi_{x}u)_{x}+\xi_{x}k_{1}u_{x} && \text{ in } Q, \label{adj}\\
{{\partial z} \over {\partial t}} + {{\partial z} \over {\partial a}}+(k_{2}(x)z_{x})_{x}-\mu_{2}(t,a, x)z =\xi h_{2}+(k_{2}\xi_{x}v)_{x}+\xi_{x}k_{2}v_{x}&& \text{ in } Q,\nonumber\\
 w(t,a, 1)=w(t,a, 0)=z(t,a, 1)=z(t,a, 0)=0 && \text{ on } (0,T)\times(0,A), \nonumber\\
\nonumber  w(T,a, x)=w_{T}(a, x) && \text{ in } Q_{A},\\
\nonumber  z(T,a, x)=z_{T}(a, x) && \text{ in } Q_{A},\\
\nonumber  w(t,A, x)=z(t,A, x)=0 && \text{ in } Q_{T}.
\end{eqnarray}
Using Proposition\ref{15} for  the inhomogeneous term $\xi(h_{1}+\mu_{3}v)+(k_{1}\xi_{x}u)_{x}+\xi_{x}k_{1}u_{x}$,
the definition of $\xi$ and Young inequality, we get the following inequality
\begin{eqnarray}\label{482}
\nonumber & & \int_{Q}(s\Theta k_{1} w_{x}^{2}+s^{3}\Theta^{3}\frac{x^{2}}{k_{1}}w^{2})e^{2s\varphi_{1}} dtdadx\\
\nonumber& &\leq  C (\int_{Q}[\xi^{2}(h_{1}+\mu_{3}v)^{2}+((k_{1}\xi_{x}u)_{x}+\xi_{x}k_{1}u_{x})^{2}]e^{2s\varphi_{1}} dtdadx
\\ \nonumber & &+ sk_{1}(1) \int_{0}^{A}\int_{0}^{T}\Theta w_{x}^{2}(t,a, 1)e^{2s\varphi_{1}(t,a, 1)} dtda)
\\ &&\leq \overline{C}\int_{Q}[\mu_{3}^{2}z^{2}+\xi^{2}h_{1}^{2}+((k_{1}\xi_{x}u)_{x}+\xi_{x}k_{1}u_{x})^{2}]e^{2s\varphi_{1}} dtdadx.
\end{eqnarray}
Thanks again to the definition of $\xi$, we have
\begin{eqnarray}\label{1003}
\nonumber& &\int_{0}^{1}((k_{1}\xi_{x}u)_{x}+\xi_{x}k_{1}u_{x})^{2}e^{2s\varphi_{1}}dx
\\ \nonumber & &\leq\int_{\omega^{'}}(8(k_{1}\xi_{x})^{2}u_{x}^{2}+2((k_{1}\xi_{x})_{x})^{2}u^{2})e^{2s\varphi_{1}}dx
\\  & & \leq C\int_{\omega^{'}}(u^{2}+u_{x}^{2})e^{2s\varphi_{1}}dx.
\end{eqnarray}
On the other hand, since $\frac{x^{2}}{k_{2}(x)}$ is non-decreasing, with the help of Hardy-Poincar\'{e} inequality
stated in \cite{Bouss}  and since $\varphi_{1}\leq\varphi_{2}$ we get
\begin{eqnarray}
\nonumber \int_{0}^{1} \mu_{3}^{2}z^{2}e^{2s\varphi_{1}}dx\leq \frac{\|\mu_{3}\|_{\infty}^{2}}{k_{2}(1)}\int_{0}^{1}\frac{k_{2}(x)}{x^{2}}(ze^{s\varphi_{2}})^{2}dx
\leq C \frac{\|\mu_{3}\|_{\infty}^{2}}{k_{2}(1)}\int_{0}^{1}k_{2}(x)((ze^{s\varphi_{2}})_{x})^{2}dx.
\end{eqnarray}
Thus, from the definition of $\psi_{2}$,  we obtain
\begin{eqnarray}
\nonumber \int_{0}^{1} \mu_{3}^{2}z^{2}e^{2s\varphi_{1}}dx\leq C \int_{0}^{1}k_{2}(x)z_{x}^{2}e^{2s\varphi_{2}}dx+C\int_{0}^{1}s^{2}\Theta^{2}\frac{x^{2}}{k_{2}(x)}z^{2}e^{2s\varphi_{2}}dx.
\end{eqnarray}
Hence, for $s$ quite large we get
\begin{eqnarray}\label{1002}
 \int_{0}^{1} \mu_{3}^{2}z^{2}e^{2s\varphi_{1}}dx\leq \frac{1}{2} \int_{0}^{1}s\Theta k_{2}(x)z_{x}^{2}e^{2s\varphi_{2}}dx+\frac{1}{2}\int_{0}^{1}s^{3}\Theta^{3}\frac{x^{2}}{k_{2}(x)}z^{2}e^{2s\varphi_{2}}dx.
\end{eqnarray}
Combining \eqref{482}, \eqref{1003} and \eqref{1002}, for $s$ quite large the following inequality holds
\begin{eqnarray}\label{carle-1}
 & & \int_{Q}(s\Theta k_{1} w_{x}^{2}+s^{3}\Theta^{3}\frac{x^{2}}{k_{1}}w^{2})e^{2s\varphi_{1}} dtdadx\\
\nonumber & & \leq \overline{C} \int_{Q}h_{1}^{2}e^{2s\varphi_{1}} dtdxda +
\frac{1}{2} \int_{Q}(s\Theta k_{2}(x) z_{x}^{2}+s^{3}\Theta^{3}\frac{x^{2}}{k_{2}(x)}z^{2})e^{2s\varphi_{2}} dtdadx
\\ \nonumber & & + C_{1}\int_{\omega^{'}}\int_{0}^{A}\int_{0}^{T}(u^{2}+u_{x}^{2})e^{2s\varphi_{1}}dtdadx.
\end{eqnarray}
 Applying the same way with $\xi h_{2}+(k_{2}\xi_{x}v)_{x}+\xi_{x}k_{2}v_{x}$ we obtain
\begin{eqnarray}\label{carle-2}
\nonumber & & \int_{Q}(s\Theta k_{2} z_{x}^{2}+s^{3}\Theta^{3}\frac{x^{2}}{k_{2}}z^{2})e^{2s\varphi_{2}} dtdadx\\
& & \leq C_{2} \int_{Q}h_{2}^{2}e^{2s\varphi_{2}} dtdxda
+C_{3}\int_{\omega^{'}}\int_{0}^{A}\int_{0}^{T}(v^{2}+v_{x}^{2})e^{2s\varphi_{2}}dtdadx.
\end{eqnarray}
Therefore, for $s$ quite large we conclude by
inequalities \eqref{carle-1} and \eqref{carle-2} and again $\varphi_{1}\leq\varphi_{2}$ that
\begin{eqnarray}\label{carle-3}
 \nonumber & & \int_{Q}(s\Theta k_{1} w_{x}^{2}+s^{3}\Theta^{3}\frac{x^{2}}{k_{1}}w^{2})e^{2s\varphi_{1}} dtdadx+
\int_{Q}(s\Theta k_{2} z_{x}^{2}+s^{3}\Theta^{3}\frac{x^{2}}{k_{2}}z^{2})e^{2s\varphi_{2}} dtdadx\\
\nonumber& & \leq C_{4} \int_{Q}(h_{1}^{2}+h_{2}^{2})e^{2s\varphi_{2}} dtdadx
+C_{5}\int_{\omega^{'}}\int_{0}^{A}\int_{0}^{T}(u^{2}+v^{2}+u_{x}^{2}+v_{x}^{2})e^{2s\varphi_{2}}dtdadx.
\end{eqnarray}
Using Caccioppoli's inequality \eqref{Caccio-1}, the last inquality becomes
\begin{eqnarray}\label{carle-4}
 \nonumber& & \int_{Q}(s\Theta k_{1} w_{x}^{2}+s^{3}\Theta^{3}\frac{x^{2}}{k_{1}}w^{2})e^{2s\varphi_{1}} dtdadx+
\int_{Q}(s\Theta k_{2} z_{x}^{2}+s^{3}\Theta^{3}\frac{x^{2}}{k_{2}}z^{2})e^{2s\varphi_{2}} dtdadx\\
& & \leq C_{6} \int_{Q}(h_{1}^{2}+h_{2}^{2})e^{2s\varphi_{2}} dtdadx
+C_{7}\int_{q}s^{2}\Theta^{2}(u^{2}+v^{2})e^{2s\varphi_{2}}dtdadx.
\end{eqnarray}
Now, let  $W:=\eta u$ and $Z:=\eta v$ with $\eta=1-\xi$. Then $W$ and $Z$ are supported in $(x_{1},1)$ and verify the following system
\begin{align}
{{\partial W} \over {\partial t}} + {{\partial W} \over {\partial a}}+(k_{1}(x)W_{x})_{x}-\mu_{1}(t,a, x)W= \mu_{3}(t,a, x)Z+\eta h_{1}+(k_{1}\eta_{x}u)_{x}+\eta_{x}k_{1}u_{x}  \quad \text{ in } Q_{x_{1}}, \label{adj-1}\\
{{\partial Z} \over {\partial t}} + {{\partial Z} \over {\partial a}}+(k_{2}(x)Z_{x})_{x}-\mu_{2}(t,a, x)Z=\eta h_{2}+(k_{2}\eta_{x}v)_{x}+\eta_{x}k_{2}v_{x} \quad \text{ in } Q_{x_{1}},\nonumber\\
W(t,a, 1)=W(t,a, x_{1})=Z(t,a, 1)=Z(t,a, x_{1})=0 \quad \text{ on } (0,T)\times(0,A), \nonumber\\
\nonumber  W(t,a, x)=W_{T}(a, x) \quad \text{ in } Q_{A},\\
\nonumber  Z(t,a, x)=Z_{T}(a, x) \quad \text{ in } Q_{A},\\
\nonumber  W(t,A, x)=Z(t,A, x)=0 \quad \text{ in } Q_{T},
\end{align}
where, $Q_{x_{1}}:=(0,T)\times(0,A)\times(x_{1},1)$. Then, the system satisfied by $W$ and $Z$ is nondegenerate.
Hence, applying Proposition \ref{Carl-nondegenerate} on the first equation of \eqref{adj-1}  for $b_{1}=x_{1}$, $b_{2}=1$ and $h:=\eta(h_{1}+\mu_{3}v)+(k_{1}\eta_{x}u)_{x}+\eta_{x}k_{1}u_{x}$, with the aid of Caccioppoli's inequality stated in \cite[Lemma 5.1]{man}, thanks to the definition of $\eta$ and Young inequality and taking $s$ quite large we obtain the following estimate
\begin{eqnarray}\label{carle-5}
\nonumber& & \int_{Q}(s^{3}\phi^{3}W^{2}+s\phi W_{x}^{2})e^{2s\Phi}dtdadx\\ \nonumber& & \leq C\left(\int_{Q}(\eta(h_{1}
+\mu_{3}v)+(k\eta_{x}u)_{x}+k\eta_{x}u_{x})^{2}e^{2s\Phi} dtdadx+ \int_{\omega}\int_{0}^{A}\int_{0}^{T} s^{3}\Theta^{3}u^{2}e^{2s\Phi} dtdadx\right)
\\ \nonumber&& \leq \widetilde{C}\left(\int_{Q}\eta^{2}(h_{1}
+\mu_{3}v)^{2}e^{2s\Phi}+((k\eta_{x}u)_{x}+k\eta_{x}u_{x})^{2}e^{2s\Phi} dtdadx+ \int_{\omega}\int_{0}^{A}\int_{0}^{T} s^{3}\Theta^{3}u^{2}e^{2s\Phi} dtdadx\right)
\\ \nonumber&& \leq \widetilde{C}(\int_{Q}\eta^{2}(h_{1}
+\mu_{3}v)^{2}e^{2s\Phi}dtdadx
+\int_{\omega^{'}}\int_{0}^{A}\int_{0}^{T}(8(k\eta_{x})^{2}u_{x}^{2}+2((k\eta_{x})_{x})^{2}u^{2})e^{2s\Phi} dtdadx\\ \nonumber&&+ \int_{\omega}\int_{0}^{A}\int_{0}^{T}s^{3}\Theta^{3} u^{2}e^{2s\Phi} dtdadx)\\ \nonumber & & \leq \widetilde{C}_{1}\left(\int_{Q}\eta^{2}(h_{1}
+\mu_{3}v)^{2}e^{2s\Phi}dtdadx
+\int_{\omega^{'}}\int_{0}^{A}\int_{0}^{T}(u_{x}^{2}+u^{2})e^{2s\Phi} dtdadx+ \int_{\omega}\int_{0}^{A}\int_{0}^{T}s^{3}\Theta^{3} u^{2}e^{2s\Phi} dtdadx\right)\\ \nonumber& & \leq \widetilde{C}_{2}\left(\int_{Q}\eta^{2}(h_{1}+\mu_{3}v)^{2}e^{2s\Phi}dtdadx+ \int_{\omega}\int_{0}^{A}\int_{0}^{T} s^{3}\Theta^{3}u^{2}e^{2s\Phi} dtdadx\right)\\ & & \leq \widetilde{C}_{3}\left(\int_{Q}(h_{1}^{2}
+\mu_{3}^{2}Z^{2})e^{2s\Phi}dtdadx+\int_{q}s^{3}\Theta^{3}u^{2}e^{2s\Phi}dtdadx\right),
\end{eqnarray}
with $\Phi$ and $\phi$ are defined in \eqref{479} and $\omega^{'}$ is defined in the beginning of the proof. On the other hand, using the fact that $x\mapsto\frac{x^{2}}{k_{2}(x)}$ is non-decreasing, Hardy-Poincar\'{e} inequality for the function $Ze^{s\Phi}$ and the definition of $\psi_{2}$ we have for $s$ quite large the following inequality
\begin{eqnarray}\label{carl}
\nonumber &&\int_{Q} \mu_{3}^{2}Z^{2}e^{2s\Phi}dx\leq c \left(\int_{Q} k_{2}(x)Z_{x}^{2}e^{2s\Phi}dtdadx+\int_{Q}s^{2}\Theta^{2}\frac{x^{2}}{k_{2}(x)}Z^{2}e^{2s\Phi}dtdadx\right)
\\ && \leq \frac{1}{2}
\int_{Q}(s^{3}\phi^{3}Z^{2}+s\phi Z_{x}^{2})e^{2s\Phi}dtdadx.
\end{eqnarray}
Therefore, injecting \eqref{carl} in \eqref{carle-5} we get
\begin{eqnarray}\label{carle-6}
& & \int_{Q}(s^{3}\phi^{3}W^{2}+s\phi W_{x}^{2})e^{2s\Phi}dtdadx\\ \nonumber & & \leq C\left(\int_{Q}h_{1}^{2}
e^{2s\Phi}dtdadx+\int_{q}s^{3}\Theta^{3}u^{2}e^{2s\Phi}dtdadx\right)+\frac{1}{2}
\int_{Q}(s^{3}\phi^{3}Z^{2}+s\phi Z_{x}^{2})e^{2s\Phi}dtdadx.
\end{eqnarray}
Replying the same argument for  the source term $h:=\eta h_{2}+(k_{2}\eta_{x}v)_{x}+\eta_{x}k_{2}v_{x}$ we infer that
\begin{eqnarray}\label{carle-7}
\int_{Q}(s^{3}\phi^{3}Z^{2}+s\phi Z_{x}^{2})e^{2s\Phi}dtdadx\leq C_{8}
\left(\int_{Q}h_{2}^{2}e^{2s\Phi}dtdadx+\int_{q}s^{3}\Theta^{3}v^{2}e^{2s\Phi}dtdadx\right).
\end{eqnarray}
Subsequently, combining \eqref{carle-6} and \eqref{carle-7} we arrive to
\begin{eqnarray}\label{carle-8}
\nonumber& &\int_{Q}[s^{3}\phi^{3}(W^{2}+Z^{2})+s\phi (W_{x}^{2}+Z_{x}^{2})]e^{2s\Phi}dtdadx\\  & &\leq C_{9}
\left(\int_{Q}(h_{1}^{2}+h_{2}^{2})e^{2s\Phi}dtdadx+\int_{q}s^{3}\Theta^{3}(u^{2}+v^{2})e^{2s\Phi}dtdadx\right).
\end{eqnarray}
Using the fact that $u=w+W$ and $v=z+Z$, $\varphi_{1}\leq\varphi_{2}\leq\Phi$, the estimates \eqref{carle-4} and \eqref{carle-8}lead to estimate \eqref{omega-estimate}.
\endproof
Using the Theorem \ref{480} for a special functions $h_{1}$ and $h_{2}$, we are ready to deduce the following result
\begin{theorem}\label{483}
Assume that the assumptions
\eqref{473} and  \eqref{474}  hold. Let $A>0$ and $T>0$ be given such that $T\in(0,\delta)$ with $\delta\in(0,A)$ small enough. Then, there exist positive constants $C$ (independent of $\delta$)  and $s_{0}$ such that
for all $s\geq s_{0}$,  every solution $(u,v)$ of \eqref{477} satisfies
\begin{eqnarray}\label{484}
\nonumber & &\int_{Q} \left(s^{3}\Theta^{3}\frac{x^{2}}{k_{1}(x)}u^{2}+s\Theta k_{1}(x)u_{x}^{2}\right)e^{2s\varphi_{1}} dtdadx
+ \int_{Q} \left(s^{3}\Theta^{3}\frac{x^{2}}{k_{2}(x)}v^{2}+s\Theta k_{2}(x)v_{x}^{2}\right)e^{2s\varphi_{2}} dtdadx
\\ & &\leq C\left(\int_{q}s^{3}\Theta^{3}(u^{2}+v^{2})e^{2s\Phi}dtdadx
+\int_{0}^{1}\int_{0}^{\delta}(u_{T}^{2}(a,x)+v_{T}^{2}(a,x))dadx\right).
\end{eqnarray}
\end{theorem}
\proof
Let $h_{1}:=-\beta_{1}(t,a,x)u(t,0,x)$ and $h_{2}:=-\beta_{2}(t,a,x)v(t,0,x)$.
\\Therefore, thanks to \eqref{omega-estimate} and \eqref{474} we have the existence of two positive constants $C$ and $s_{0}$ such that, for all $s\geq s_{0}$, the following inequality holds
\begin{eqnarray}
\nonumber
& &
s^{3}\int_{Q}  \Theta^{3}\left(\frac{x^{2}}{k_{1}(x)}u^{2}e^{2s\varphi_{1}}+\frac{x^{2}}{k_{2}(x)}v^{2}e^{2s\varphi_{2}}\right) dtdadx + s\int_{Q}  \Theta \left(k_{1}(x)u_{x}^{2}e^{2s\varphi_{1}}+k_{2}(x)v_{x}^{2}e^{2s\varphi_{2}}\right) dtdadx \\ \nonumber& & \leq C
\left(\int_{Q}((\beta_{1})^{2} u^{2}(t,0,x)+(\beta_{2})^{2} v^{2}(t,0,x))e^{2s\Phi}dtdadx+\int_{q}s^{3}\Theta^{3}(u^{2}+v^{2})e^{2s\Phi}dtdadx\right) \label{betas2}
\\ && \leq\widetilde{C}_{1}
\left(\int_{0}^{1}\int_{0}^{T}( u^{2}(t,0,x)+v^{2}(t,0,x))dtdadx+\int_{q}s^{3}\Theta^{3}(u^{2}+v^{2})e^{2s\Phi}dtdadx\right)
\end{eqnarray}
Set $U(t,a,x)=u(T-t,A-a,x)$ and $V(t,a,x)=v(T-t,A-a,x)$. Then, one has
\begin{align}
&{{\partial U} \over {\partial t}} + {{\partial U} \over {\partial a}}-(k_{1}(x)U_{x})_{x}+\mu_{1}(T-t,A-a, x)U+\mu_{3}(T-t,A-a, x)V =\beta_{1}(T-t,A-a,x)U(t,A,x),\nonumber\\
& U(t,a, 1)=U(t,a, 0)=0, \label{302}\\
&\nonumber  U(0,a, x)=U_{0}(a,x)=u_{T}(A-a, x),\\
&\nonumber  U(t,0,x)=0,
\end{align}
where $V$ is the solution of
\begin{align}
&{{\partial V} \over {\partial t}} + {{\partial V} \over {\partial a}}-(k_{2}(x)V_{x})_{x}+\mu_{2}(T-t,A-a, x)V =\beta_{2}(T-t,A-a,x)V(t,A,x),\nonumber\\
& V(t,a, 1)=V(t,a, 0)=0, \label{303}\\
&\nonumber  V(0,a, x)=V_{0}(a,x)=v_{T}(A-a, x),\\
&\nonumber  V(t,0,x)=0.
\end{align}
Integrating along the characteristic lines, we get respectively the implicit formulas for the solutions $U$ of \eqref{302}  and $V$ of \eqref{303} given by
\begin{equation}\label{line-1}
\begin{cases}
U(t,a,\cdot)=\int_{0}^{a}S(a-l)(\beta_{1}(T-t,A-l,\cdot)U(t,A,\cdot)-\mu_{3}(T-t,A-l,\cdot)V(t,l,\cdot))dl, \\ \text{ if } t>a\\
U(t,a,\cdot)=S(t)U_{0}(a-t,\cdot)+\int_{0}^{t}S(t-l)(\beta_{1}(T-l,A-a,\cdot)U(l,A,\cdot)-\mu_{3}(T-l,A-a,\cdot)V(l,a,\cdot))dl, \\ \text{ if } t\leq a,
\end{cases}
\end{equation}
and
\begin{equation}\label{453}
\begin{cases}
V(t,a,\cdot)=\int_{0}^{a}\mathrm{L}(a-l)\beta_{2}(T-t,A-l,\cdot)V(t,A,\cdot)dl, & \text{ if } t>a\\
V(t,a,\cdot)=\mathrm{L}(t)V_{0}(a-t,\cdot)+\int_{0}^{t}\mathrm{L}(t-l)\beta_{2}(T-l,A-a,\cdot)V(l,A,\cdot)dl, & \text{ if } t\leq a,
\end{cases}
\end{equation}
where $(S(t))_{t\geq 0}$ and $(\mathrm{L}(t))_{t\geq 0}$ are  the bounded semigroups  generated respectively by the operators $A_{4}U=-(k_{1}U_{x})_{x}+\mu_{1}(T-t,A-a, x)U$ and $A_{7}V=-(k_{2}V_{x})_{x}+\mu_{2}(T-t,A-a, x)V$.\\
Hence, after a careful computations, \eqref{line-1} and \eqref{453} become respectively
\begin{equation}\label{line-3}
\begin{cases}
u(t,a,\cdot)=\int_{0}^{A-a}S(A-a-l)(\beta_{1}(t,A-l,\cdot)u(t,0,\cdot)-\mu_{3}(t,A-l,\cdot)v(t,A-l,\cdot))dl, \\ \text{ if } a>t+(A-T)\\
u(t,a,\cdot)=S(T-t)u_{T}(T+(a-t),\cdot)+\int_{t}^{T}S(l-t)(\beta_{1}(l,a,\cdot)u(l,0,\cdot)-\mu_{3}(l,a,\cdot)v(l,a,\cdot))dl, \\ \text{ if } a\leq t+(A-T),
\end{cases}
\end{equation}
\begin{equation}\label{line-2}
\begin{cases}
v(t,a,\cdot)=\int_{0}^{A-a}\mathrm{L}(A-a-l)\beta_{2}(t,A-l,\cdot)v(t,0,\cdot)dl, & \text{ if } a>t+(A-T)\\
v(t,a,\cdot)=\mathrm{L}(T-t)v_{T}(T+(a-t),\cdot)+\int_{t}^{T}\mathrm{L}(l-t)\beta_{2}(l,a,\cdot)v(l,0,\cdot)dl, & \text{ if } a\leq t+(A-T),
\end{cases}
\end{equation}
Thus, by the third hypothesis in \eqref{474} on $\beta_{1}$ and $\beta_{2}$ one has
\begin{equation}\label{line-4}
\begin{cases}
u(t,0,\cdot)=S(T-t)u_{T}(T-t,\cdot)-\int_{t}^{T}S(l-t)\mu_{3}(l,0,\cdot)v(l,0,\cdot)dl,\\
v(t,0,\cdot)=\mathrm{L}(T-t)v_{T}(T+(a-t),\cdot).
\end{cases}
\end{equation}

Subsequently, by \eqref{betas2} we deduce that
\begin{eqnarray}\label{betas3}
\nonumber && s^{3}\int_{Q}  \Theta^{3}\left(\frac{x^{2}}{k_{1}(x)}u^{2}e^{2s\varphi_{1}}+\frac{x^{2}}{k_{2}(x)}v^{2}e^{2s\varphi_{2}}\right) dtdadx
+ s\int_{Q}  \Theta \left(k_{1}(x)u_{x}^{2}e^{2s\varphi_{1}}+k_{2}(x)v_{x}^{2}e^{2s\varphi_{2}}\right) dtdadx
\\ && \leq \widehat{C}_{1}\left(\int_{q}s^{3}\Theta^{3}(u^{2}+v^{2})e^{2s\Phi}dtdadx
+\int_{0}^{1}\int_{0}^{\delta}(u_{T}^{2}(a,x)+v_{T}^{2}(a,x))dadx\right),
\end{eqnarray}
since $(S(t))_{t\geq0}$ and $(\mathrm{L}(t))_{t\geq0}$ are a bounded semigroups, $\mu_{3}\in L^{\infty}(Q)$ and $T\in(0, \delta)$.\\
Then the thesis follows.
\endproof
We come now to the more challenging point and the novelty of this contribution which is the following $\omega$-Carleman type inequality. Such an estimate plays a crucial role to obtain the null controllability of population dynamics cascade system with one control force.
\begin{theorem}\label{theo-1}
Let
\eqref{473} and  \eqref{474} be verified. Let $A>0$ and $T>0$ be given such that $T\in(0,\delta)$ with $\delta\in(0,A)$ small enough. Assume that there exists  a positive constant $\nu$ such that
\begin{equation}\label{mu-3}
\mu_{3}\geq \nu \hspace{0.5cm}\text{ on } [0,T]\times[0,A]\times\omega_{1} \hspace{0.25cm}\text{ for some }  \omega_{1}\Subset \omega,
\end{equation}

Then every solution $(u,v)$ of \eqref{477} satisfies
\begin{align}\label{Car}
\nonumber& &\int_{Q} \left(s^{3}\Theta^{3}\frac{x^{2}}{k_{1}(x)}u^{2}+s\Theta k_{1}(x)u_{x}^{2}\right)e^{2s\varphi_{1}} dtdadx
+ \int_{Q} \left(s^{3}\Theta^{3}\frac{x^{2}}{k_{2}(x)}v^{2}+s\Theta k_{2}(x)v_{x}^{2}\right)e^{2s\varphi_{2}} dtdadx
\\ & &\leq C_{\delta}\left( \int_{q}u^{2}dtdadx
+\int_{0}^{1}\int_{0}^{\delta}(u_{T}^{2}(a,x)+v_{T}^{2}(a,x))dadx\right).
\end{align}
\end{theorem}
This inequality is an immediate outcome of Theorem \ref{483} applied to $\omega_{1}$ and the following lemma
(see for instance \cite{Haj1} and the references therein).
\begin{lemma}\label{lemma-1}
Assume that
\eqref{473} and  \eqref{474}  hold and let $A>0$ and $T>0$ be given such that $T\in(0,\delta)$ with $\delta\in(0,A)$ small enough. we suppose also that \eqref{mu-3} holds. Then, for all $\epsilon>0$
there exist two positive constants $C$ and $M_{\epsilon}$ such that for every solution $(u,v)$ of \eqref{477} the following inequality
is satisfied
\begin{eqnarray}\label{final}
\nonumber &&\int_{\omega_{1}}\int_{0}^{A}\int_{0}^{T} s^{3}\Theta^{3}v^{2}e^{2s\Phi}dtdadx\leq\epsilon C\left(\int_{Q}s^{3}\Theta^{3}\frac{x^{2}}{k_{2}}v^{2}e^{2s\varphi_{2}}dtdadx
+\int_{Q}s\Theta k_{2}(x)v_{x}^{2}e^{2s\varphi_{2}}dtdadx\right)
\\ &&+M_{\epsilon}\left(\int_{\omega}\int_{0}^{A}\int_{0}^{T}u^{2}dtdadx
+\int_{0}^{1}\int_{0}^{\delta}(u_{T}^{2}(a,x)+v_{T}^{2}(a,x))dadx\right).
\end{eqnarray}
\end{lemma}
\proof
Let $\chi:\mathbb{R}\rightarrow \mathbb{R}$ be the non-negative cut-off function defined as follows
\begin{equation}\label{cut-off func 1}
\begin{cases}
\chi\in\mathcal{C}^{\infty}(0,1),\\
supp(\chi)\subset\omega,\\
\chi\equiv 1 \hspace{0.5cm} \text{ on } \omega_{1}.
\end{cases}
\end{equation}
Recall that $\omega=(x_{1},x_{2})$. Multiplying the first equation of \eqref{477} by $\chi s^{3}\Theta^{3}ve^{2s\Phi}$ and after an integration by parts, we get
\begin{eqnarray}
\nonumber&& \int_{Q}\chi s^{3}\Theta^{3}ve^{2s\Phi}u_{t}dtdadx=-\int_{Q}(3+2s\Phi)\chi s^{3}\Theta_{t}\Theta^{2}uve^{2s\Phi}dtdadx
-\int_{Q}\chi s^{3}\Theta^{3}uv_{t}e^{2s\Phi}dtdadx.\\ \nonumber&& \int_{Q}\chi s^{3}\Theta^{3}ve^{2s\Phi}u_{a}dtdadx=-\int_{Q}(3+2s\Phi)\chi s^{3}\Theta_{a}\Theta^{2}uve^{2s\Phi}dtdadx
-\int_{Q}\chi s^{3}\Theta^{3}uv_{a}e^{2s\Phi}dtdadx.\\ \nonumber&& \int_{Q}\chi s^{3}\Theta^{3}ve^{2s\Phi}(k_{1}u_{x})_{x}dtdadx=
-\int_{Q}\chi s^{3}\Theta^{3}k_{1}e^{2s\Phi}u_{x}v_{x}dtdadx+\int_{Q} s^{3}\Theta^{3}k_{1}(\chi e^{2s\Phi})_{x}uv_{x}dtdadx
\\ \nonumber&&+\int_{Q}s^{3}\Theta^{3}(k_{1}(\chi e^{2s\Phi})_{x})_{x}uv dtdadx.
\\ \nonumber&& -\int_{Q}\chi s^{3}\Theta^{3}ve^{2s\Phi}\mu_{1}u dtdadx=-\int_{Q}\chi s^{3}\Theta^{3}\mu_{1}uve^{2s\Phi}dtdadx.
\\ \nonumber &&-\int_{Q}\chi s^{3}\Theta^{3}ve^{2s\Phi}\mu_{3}v dtdadx=-\int_{Q}\chi s^{3}\Theta^{3}\mu_{3}v^{2}e^{2s\Phi}dtdadx.
\end{eqnarray}
Then, summing all these identities side by side, using the second equation of \eqref{477} and integrating again by parts
 \begin{eqnarray}\label{int-by-parts}
& & \int_{Q}\chi s^{3}\Theta^{3}\mu _{3}v^{2}e^{2s\Phi}dtdadx=I_{1}+I_{2}+I_{3}+I_{4}+I_{5},
\end{eqnarray}
where,
$I_{1}:=\int_{Q}\chi s^{3}\Theta^{3}\beta_{1}v u(t,0,x) e^{2s\Phi}dtdadx$,\\
$I_{2}:=-\int_{Q} \left((3+2s\Phi) s^{3}\Theta_{t}\Theta^{2}+(3+2s\Phi) s^{3}\Theta_{a}\Theta^{2}+\mu_{1}s^{3}\Theta^{3}+\mu_{2}s^{3}\Theta^{3}\right)\chi e^{2s\Phi}uv dtdadx
\\+\int_{Q}s^{3}\Theta^{3}(k_{1}(\chi e^{2s\Phi})_{x})_{x}uvdtdadx$,\\
$I_{3}:=\int_{Q} \chi s^{3}\Theta^{3} \beta_{2}uv(t,0,x)e^{2s\Phi}dtdadx$,
$I_{4}:=\int_{Q}s^{3}\Theta^{3}(k_{1}-k_{2})(x)uv_{x}(\chi e^{2s\Phi})_{x}dtdadx$,\\
$I_{5}:=-\int_{Q}\chi s^{3}\Theta^{3}(k_{1}+k_{2})(x)u_{x}v_{x} e^{2s\Phi}dtdadx$.\\
On one hand, we have by Young inequality and definition of $\chi$
\begin{eqnarray}\label{I-5}
\nonumber & &I_{5}\leq \epsilon\int_{Q}s\Theta k_{2}(x)v_{x}^{2}e^{2s\varphi_{2}}dtdadx
+\frac{1}{4\epsilon}\int_{Q}\frac{\chi^{2}s^{5}\Theta^{5}(k_{1}+k_{2})^{2}u_{x}^{2}e^{2s(2\Phi-\varphi_{2})}}{k_{2}}dtdadx
\\ \nonumber && \leq \epsilon\int_{Q}s\Theta k_{2}(x)v_{x}^{2}e^{2s\varphi_{2}}dtdadx
\\ &&+\frac{\max_{[0,1]}(k_{1}+k_{2})^{2}}{4\epsilon\min_{\omega}k_{2}}\int_{Q}\chi s^{5}\Theta^{5}u_{x}^{2}e^{2s(2\Phi-\varphi_{2})}dtdadx.
\end{eqnarray}
Put $L:=\int_{Q}\chi s^{5}\Theta^{5}u_{x}^{2}e^{2s(2\Phi-\varphi_{2})}dtdadx$. To increase $I_{5}$, we will find an upper
bound of $L$. To do this, we multiply the first equation of \eqref{477} by
$\frac{\chi s^{5}\Theta^{5}e^{2s(2\Phi-\varphi_{2})}}{k_{1}}u$ and after integration by parts
\begin{eqnarray}
\nonumber && \int_{Q}\frac{\chi s^{5}\Theta^{5}e^{2s(2\Phi-\varphi_{2})}}{k_{1}}uu_{t}dtdadx=
-\frac{1}{2}\int_{Q}\frac{s^{5}\chi}{k_{1}}\Theta^{4}\Theta_{t}(5+2s(2\Phi-\varphi_{2}))e^{2s(2\Phi-\varphi_{2})}u^{2}dtdadx.
\\ \nonumber && \int_{Q}\frac{\chi s^{5}\Theta^{5}e^{2s(2\Phi-\varphi_{2})}}{k_{1}}uu_{a}dtdadx=
-\frac{1}{2}\int_{Q}\frac{s^{5}\chi}{k_{1}}\Theta^{4}\Theta_{a}(5+2s(2\Phi-\varphi_{2}))e^{2s(2\Phi-\varphi_{2})}u^{2}dtdadx.
\\ \nonumber && \int_{Q}\frac{\chi s^{5}\Theta^{5}e^{2s(2\Phi-\varphi_{2})}}{k_{1}}u(k_{1}u_{x})_{x}dtdadx=
-\int_{Q}\chi s^{5}\Theta^{5}u_{x}^{2}e^{2s(2\Phi-\varphi_{2})}dtdadx
\\ \nonumber &&+\frac{1}{2}\int_{Q}s^{5}\Theta^{5}\left(k_{1}\left(\frac{\chi e^{2s(2\Phi-\varphi_{2})}}{k_{1}}\right)_{x}\right)_{x}u^{2}dtdadx.
\\ \nonumber && -\int_{Q}\frac{\chi s^{5}\Theta^{5}e^{2s(2\Phi-\varphi_{2})}}{k_{1}}u\mu_{1}udtdadx=
-\int_{Q}\frac{\chi s^{5}\Theta^{5}e^{2s(2\Phi-\varphi_{2})}}{k_{1}}\mu_{1}u^{2}dtdadx.
\\ \nonumber && -\int_{Q}\frac{\chi s^{5}\Theta^{5}e^{2s(2\Phi-\varphi_{2})}}{k_{1}}u\mu_{3}vdtdadx=
-\int_{Q}\frac{\chi s^{5}\Theta^{5}e^{2s(2\Phi-\varphi_{2})}}{k_{1}}\mu_{3}uv dtdadx.
\end{eqnarray}
Hence, adding these equalities side by side we get
\begin{eqnarray}\label{L}
L=L_{1}+L_{2}+L_{3},
\end{eqnarray}
where,
$L_{1}:=\int_{Q}\frac{\chi s^{5}\Theta^{5}e^{2s(2\Phi-\varphi_{2})}}{k_{1}}\beta_{1}uu(t,0,x)dtdadx.$\\
$L_{2}:=-\int_{Q}\frac{\chi s^{5}\Theta^{5}e^{2s(2\Phi-\varphi_{2})}}{k_{1}}\mu_{3}uv dtdadx.$\\
$L_{3}:=-\int_{Q}\left(\frac{\chi s^{5}\Theta^{5}}{k_{1}}\mu_{1}
+\frac{1}{2}\frac{s^{5}\chi}{k_{1}}\Theta^{5}\Theta_{t}(\frac{5}{\Theta}+2s(2\Psi-\psi_{2}))
+\frac{1}{2}\frac{s^{5}\chi}{k_{1}}\Theta^{5}\Theta_{a}(\frac{5}{\Theta}+2s(2\Psi-\psi_{2}))\right)e^{2s(2\Phi-\varphi_{2})}u^{2}dtdadx
\\+\frac{1}{2}\int_{Q}s^{5}\Theta^{5}\left(k_{1}\left(\frac{\chi e^{2s(2\Phi-\varphi_{2})}}{k_{1}}\right)_{x}\right)_{x}u^{2}dtdadx.$
\\ The assumptions in \eqref{474} on $\beta_{1}$ together with Young inequality, Lemma \ref{lemma-4.3},
the definitions of $\chi$ and $\Theta$, the fact that the function $x\mapsto\frac{k_{2}}{x^{2}}$ is non-increasing,
$|\Theta_{t}|\leq C \Theta^{2}$  and $|\Theta_{a}|\leq \widetilde{C} \Theta^{2}$ and
\begin{eqnarray}\label{sup}
\sup_{(t,a,x)\in Q}s^{p}\Theta^{p}e^{2s(2\Phi-\varphi_{2})}<+\infty \hspace{0.25cm} \text{ for } p\in\mathbb{R},
\end{eqnarray}
 lead to
\begin{align}\label{L-1}
\nonumber && L_{1}\leq\frac{1}{4\epsilon}\int_{Q}\frac{\chi s^{5}\Theta^{5}e^{2s(2\Phi-\varphi_{2})}}{(k_{1})^{2}}u^{2}dtdadx
+\epsilon\int_{Q}\chi s^{5}\Theta^{5}e^{2s(2\Phi-\varphi_{2})}(\beta_{1})^{2}u^{2}(t,0,x)dtdadx
\\ \nonumber && \leq \frac{\widetilde{K}_{1}}{4\epsilon}\int_{Q}\chi s^{5}\Theta^{5}e^{2s(2\Phi-\varphi_{2})}u^{2}dtdadx+
\epsilon K_{1}\int_{0}^{1}\int_{0}^{A}\int_{T-\delta}^{T}\chi u^{2}(t,0,x)dtdadx
\\ && \leq\frac{\widetilde{K}_{1}}{4\epsilon}\int_{Q}\chi s^{5}\Theta^{5}e^{2s(2\Phi-\varphi_{2})}u^{2}dtdadx+
\epsilon K_{2}
\int_{0}^{1}\int_{0}^{\delta}\chi u_{T}^{2}(a,x)dadx
\end{align}
and
\begin{eqnarray}\label{L-2}
\nonumber && L_{2}\leq\epsilon^{2}\int_{Q}\frac{x^{2}}{k_{2}}s^{3}\Theta^{3}e^{2s\varphi_{2}}v^{2}dtdadx
+\frac{1}{4\epsilon^{2}}\int_{Q}\chi^{2}\frac{s^{7}\Theta^{7}}{(k_{1})^{2}}e^{2s(4\Phi-3\varphi_{2})}\frac{k_{2}}{x^{2}}(\mu_{3})^{2}u^{2}dtdadx
\\&& \leq\epsilon^{2}\int_{Q}\frac{x^{2}}{k_{2}}s^{3}\Theta^{3}e^{2s\varphi_{2}}v^{2}dtdadx
+\frac{K_{4}}{4\epsilon^{2}}\int_{\omega}\int_{0}^{A}\int_{0}^{T} s^{7}\Theta^{7}e^{2s(4\Phi-3\varphi_{2})}u^{2}dtdadx,
\end{eqnarray}
and
\begin{eqnarray}\label{L-3}
 |L_{3}|\leq K_{5}\int_{\omega}\int_{0}^{A}\int_{0}^{T} s^{7}\Theta^{7}e^{2s(2\Phi-\varphi_{2})}u^{2}dtdadx,
\end{eqnarray}
where $K_{4}=\frac{\|\mu_{3}\|_{\infty}^{2}k_{2}(x_{1})}{(x_{1})^{2}\min_{\omega}k_{1}}$. On the other hand, by Lemma \ref{interval-2} we have
\begin{eqnarray}\label{equivalence}
e^{2s(2\Phi-\varphi_{2})}\leq e^{2s(4\Phi-3\varphi_{2})}.
\end{eqnarray}
Then, combining relations \eqref{L}, \eqref{L-1}, \eqref{L-2} and \eqref{L-3} we conclude
\begin{eqnarray}\label{L-4}
\nonumber && L\leq\epsilon^{2}\int_{Q}\frac{x^{2}}{k_{2}}s^{3}\Theta^{3}e^{2s\varphi_{2}}v^{2}dtdadx
+K_{\epsilon}\int_{\omega}\int_{0}^{A}\int_{0}^{T} s^{7}\Theta^{7}e^{2s(4\Phi-3\varphi_{2})}u^{2}dtdadx
\\ &&+ \epsilon K_{2}\int_{0}^{1}\int_{0}^{\delta} u_{T}^{2}(a,x)dadx.
\end{eqnarray}
Hence, by \eqref{I-5} and \eqref{L-4} we deduce
\begin{eqnarray}\label{L-5}
\nonumber &&I_{5}\leq\epsilon C\left(\int_{Q}\frac{x^{2}}{k_{2}}s^{3}\Theta^{3}e^{2s\varphi_{2}}v^{2}dtdadx
+\int_{Q}s\Theta k_{2}(x)v_{x}^{2}e^{2s\varphi_{2}}dtdadx\right)
\\ &&+K_{\epsilon}^{1}\int_{\omega}\int_{0}^{A}\int_{0}^{T} s^{7}\Theta^{7}e^{2s(4\Phi-3\varphi_{2})}u^{2}dtdadx
+K_{2}\int_{0}^{1}\int_{0}^{\delta} u_{T}^{2}(a,x)dadx.
\end{eqnarray}
where $K_{\epsilon}^{1}$ is a positive constants that depend on $\epsilon$. Similarly, we will find an upper bounds of $I_{1}$, $I_{2}$, $I_{3}$ and $I_{4}$.
Firstly, we will start by $I_{2}$. One has the following relations
\begin{eqnarray}\label{I-2-1}
\nonumber &&\left|\int_{Q}\chi(3+2s\Phi)s^{3}\Theta_{t}\Theta^{2}e^{2s\Phi}uv dtdadx\right|
\leq\int_{Q}\chi|3+2s\Phi|s^{3}|\Theta_{t}|\Theta^{2}e^{2s\Phi}|uv|dtdadx
\\ \nonumber &&\leq C\int_{Q}\chi|3+2s\Phi|s^{3}\Theta^{4}e^{2s\Phi}|uv|dtdadx
\\ && \leq \epsilon\int_{Q}s^{3}\Theta^{3}\frac{x^{2}}{k_{2}}e^{2s\varphi_{2}}v^{2}dtdadx
+C_{\epsilon}\int_{\omega}\int_{0}^{A}\int_{0}^{T}s^{5}\Theta^{5}e^{2s(2\Phi-\varphi_{2})}u^{2}dtdadx,
\end{eqnarray}

\begin{eqnarray}\label{I-2-2}
\nonumber &&\left|\int_{Q}\chi(3+2s\Phi)s^{3}\Theta_{a}\Theta^{2}e^{2s\Phi}uv dtdadx\right|
\\ && \leq \epsilon\int_{Q}s^{3}\Theta^{3}\frac{x^{2}}{k_{2}}e^{2s\varphi_{2}}v^{2}dtdadx
+C_{\epsilon}^{1}\int_{\omega}\int_{0}^{A}\int_{0}^{T}s^{5}\Theta^{5}e^{2s(2\Phi-\varphi_{2})}u^{2}dtdadx,
\end{eqnarray}

\begin{eqnarray}\label{I-2-3}
\nonumber &&\left|\int_{Q}\chi(\mu_{1}+\mu_{2})s^{3}\Theta^{3}e^{2s\Phi}uv dtdadx\right|
\\ && \leq \epsilon\int_{Q}s^{3}\Theta^{3}\frac{x^{2}}{k_{2}}e^{2s\varphi_{2}}v^{2}dtdadx
+C_{\epsilon}^{2}\int_{\omega}\int_{0}^{A}\int_{0}^{T}s^{3}\Theta^{3}e^{2s(2\Phi-\varphi_{2})}u^{2}dtdadx,
\end{eqnarray}
\begin{eqnarray}\label{I-2-4}
\nonumber &&\left|\int_{Q}s^{3}\Theta^{3}(k_{1}(\chi e^{2s\Phi})_{x})_{x}uv dtdadx\right|
\\ \nonumber && \leq \epsilon\int_{Q}s^{3}\Theta^{3}\frac{x^{2}}{k_{2}}e^{2s\varphi_{2}}v^{2}dtdadx
+\frac{1}{4\epsilon}\int_{Q}s^{3}\Theta^{3}\frac{k_{2}}{x^{2}}(k_{1}(\chi e^{2s\Phi})_{x})_{x}^{2}e^{-2s\varphi_{2}}u^{2}dtdadx
\\ \nonumber && \leq \epsilon\int_{Q}s^{3}\Theta^{3}\frac{x^{2}}{k_{2}}e^{2s\varphi_{2}}v^{2}dtdadx
+\frac{C_{2}}{4\epsilon}\int_{Q}s^{3}\Theta^{3}\frac{k_{2}}{x^{2}}(\chi^{2}+\chi_{x}^{2}+\chi_{xx}^{2})e^{2s(2\Phi-\varphi_{2})}u^{2}dtdadx
\\ && \leq \epsilon\int_{Q}s^{3}\Theta^{3}\frac{x^{2}}{k_{2}}e^{2s\varphi_{2}}v^{2}dtdadx
+C_{\epsilon}^{3}\int_{\omega}\int_{0}^{A}\int_{0}^{T}s^{3}\Theta^{3}e^{2s(2\Phi-\varphi_{2})}u^{2}dtdadx,
\end{eqnarray}
Hence, summing inequalities \eqref{I-2-1}, \eqref{I-2-2}, \eqref{I-2-3} and \eqref{I-2-4} we obtain
\begin{eqnarray}\label{I-2}
 I_{2}\leq 4\epsilon\int_{Q}s^{3}\Theta^{3}\frac{x^{2}}{k_{2}}e^{2s\varphi_{2}}v^{2}dtdadx
+C_{\epsilon}^{4}\int_{\omega}\int_{0}^{A}\int_{0}^{T}s^{5}\Theta^{5}e^{2s(2\Phi-\varphi_{2})}u^{2}dtdadx.
\end{eqnarray}
For the rest of integrals,
\begin{align}\label{I-1}
\nonumber I_{1}=\int_{Q}\chi s^{3}\Theta^{3}\beta_{1}v u(t,0,x) e^{2s\Phi}dtdadx
\\  \leq\epsilon\int_{Q}s^{3}\Theta^{3}\frac{x^{2}}{k_{2}}e^{2s\varphi_{2}}v^{2}dtdadx
+C_{\epsilon}^{5}\int_{0}^{1}\int_{0}^{\delta} u_{T}^{2}(a,x)dadx.
\end{align}
\begin{align}\label{I-3}
\nonumber I_{3}=\int_{Q}\chi s^{3}\Theta^{3}\beta_{2}u v(t,0,x) e^{2s\Phi}dtdadx
\\  \leq\epsilon\int_{0}^{1}\int_{0}^{\delta}v_{T}^{2}(a,x)dadx
+\frac{1}{4\epsilon}\int_{\omega}\int_{0}^{A}\int_{0}^{T}s^{7}\Theta^{7}e^{2s(2\Phi-\varphi_{2})}u^{2}dtdadx.
\end{align}
\begin{eqnarray}\label{I-4}
\nonumber I_{4}&=&\int_{Q}s^{3}\Theta^{3}(k_{1}-k_{2})(x)uv_{x}(\chi e^{2s\Phi})_{x}dtdadx
\\ \nonumber &=&\int_{Q}s^{3}\Theta^{3}(k_{1}-k_{2})(x)uv_{x}(\chi_{x}+2s\Phi_{x}\chi)e^{2s\Phi}dtdadx
\\ \nonumber &&\leq\epsilon\int_{Q}s\Theta k_{2}v_{x}^{2}e^{2s\varphi_{2}}dadx+
\frac{1}{4\epsilon}\int_{Q}s^{5}\Theta^{5}\frac{(k_{1}-k_{2})^{2}}{k_{2}}(\chi_{x}+2s\Phi_{x}\chi)^{2}e^{2s(2\Phi-\varphi_{2})}u^{2}dtdadx
\\ && \leq\epsilon\int_{Q}s\Theta k_{2}v_{x}^{2}e^{2s\varphi_{2}}dadx
+C_{\epsilon}^{6}\int_{\omega}\int_{0}^{A}\int_{0}^{T}s^{7}\Theta^{7}e^{2s(2\Phi-\varphi_{2})}u^{2}dtdadx.
\end{eqnarray}
Subsequently, combining \eqref{L-5}, \eqref{I-2}, \eqref{I-1}, \eqref{I-3}, \eqref{I-4} and using again \eqref{equivalence}
\begin{eqnarray}\label{final-1}
\nonumber &&\int_{Q}\chi s^{3}\Theta^{3}\mu _{3}v^{2}e^{2s\Phi}dtdadx\leq\epsilon C_{7}\left(\int_{Q}s^{3}\Theta^{3}\frac{x^{2}}{k_{2}}v^{2}e^{2s\varphi_{2}}dtdadx
+\int_{Q}s\Theta k_{2}(x)v_{x}^{2}e^{2s\varphi_{2}}dtdadx\right)
\\ \nonumber&&+C_{\epsilon}^{8}\int_{\omega}\int_{0}^{A}\int_{0}^{T}s^{7}\Theta^{7}e^{2s(4\Phi-3\varphi_{2})}u^{2}dtdadx
+C_{\epsilon}^{9}\int_{0}^{1}\int_{0}^{\delta}(u_{T}^{2}(a,x)+v_{T}^{2}(a,x))dadx.
\end{eqnarray}
Finally, the hypothesis \eqref{mu-3}, the definition of $\chi$ and the relation
\begin{eqnarray}\label{sup-1}
\sup_{(t,a,x)\in Q}s^{p}\Theta^{p}e^{2s(4\Phi-3\varphi_{2})}<+\infty \hspace{0.25cm} \text{ for } p\in\mathbb{R},
\end{eqnarray}
 yield
\begin{eqnarray}\label{final}
\nonumber &&\int_{\omega_{1}}\int_{0}^{A}\int_{0}^{T}s^{3}\Theta^{3}v^{2}e^{2s\Phi}dtdadx\leq\epsilon C_{10}\left(\int_{Q}s^{3}\Theta^{3}\frac{x^{2}}{k_{2}}v^{2}e^{2s\varphi_{2}}dtdadx
+\int_{Q}s\Theta k_{2}(x)v_{x}^{2}e^{2s\varphi_{2}}dtdadx\right)
\\&&+C_{\epsilon}^{11}\left(\int_{\omega}\int_{0}^{A}\int_{0}^{T}u^{2}dtdadx+\int_{0}^{1}\int_{0}^{\delta}(u_{T}^{2}(a,x)+v_{T}^{2}(a,x))dadx\right),
\end{eqnarray}
which finishes the proof.
\endproof
The above Carleman estimate can be used in a standard way to obtain the null controllability
of the cascade system with one control force. This will be reached showing an observability inequality of the adjoint system.

\section{Observability inequality and  null controllability results}\label{third-section-5}
This paragraph is devoted to the observability inequality of system \eqref{477} and then the null
controllability result of system \eqref{470}. We start to show our observability inequality whose proof is based
essentially on Carleman estimate \eqref{Car} and Hardy-Poincar\'{e} inequality.
\begin{proposition}\label{485}
Assume that \eqref{473} and \eqref{474} hold. Suppose also that \eqref{mu-3} is fulfilled and let $A>0$ and $T>0$ be given such that $T\in(0,\delta)$ with $\delta\in(0,A)$ small enough. Then, there exists a positive constant $C_{\delta}$ such that  for every solution   $(u,v)$ of \eqref{477},
the  following observability inequality is satisfied
\begin{align}\label{486}
 \int_{0}^{1}\int_{0}^{A} (u^{2}(0,a,x)+v^{2}(0,a,x))dadx \leq C_{\delta} \left(\int_{q}u^{2}dtdadx+\int_{0}^{1}\int_{0}^{\delta}(u_{T}^{2}(a,x)+v^{2}_{T}(a,x))dadx\right).
\end{align}
\end{proposition}
\proof
Then for $\kappa>0$ to be defined later, $\widetilde{u}=e^{\kappa t}u$  and $\widetilde{v}=e^{\kappa t}v$ are respectively a solutions of
\begin{eqnarray}\label{obser-1}
{{\partial \widetilde{u}} \over {\partial t}} + {{\partial \widetilde{u}} \over {\partial a}}+(k_{1}(x)\widetilde{u}_{x})_{x}-\mu_{1}(t,a, x)\widetilde{u}=\mu_{3}(t,a, x)\widetilde{v}-\beta_{1}\widetilde{u}(t,0,x)  && \text{ in } Q, \label{inter-1}\\
\widetilde{u}(t,a, 1)=\widetilde{u}(t,a, 0)=0 && \text{ on } (0,T)\times(0,A), \nonumber \\
\nonumber  \widetilde{u}(T,a, x)=e^{\kappa T}  u_{T}(a, x) && \text{ in } Q_{A},\\
\nonumber  \widetilde{u}(t,A, x)=0 && \text{ in } Q_{T},
\end{eqnarray}
and
\begin{eqnarray}\label{obser-2}
{{\partial \widetilde{v}} \over {\partial t}} + {{\partial \widetilde{v}} \over {\partial a}}+(k_{2}(x)\widetilde{v}_{x})_{x}-\mu_{2}(t,a, x)\widetilde{v}=-\beta_{2}\widetilde{v}(t,0,x)  && \text{ in } Q, \label{inter-2}\\
\widetilde{v}(t,a, 1)=\widetilde{v}(t,a, 0)=0 && \text{ on } (0,T)\times(0,A), \nonumber \\
\nonumber  \widetilde{v}(T,a, x)=e^{\kappa T}v_{T}(a, x) && \text{ in } Q_{A},\\
\nonumber  \widetilde{v}(t,A, x)=0 && \text{ in } Q_{T},
\end{eqnarray}
where, $u$ and $v$ are respectively the solutions of
\begin{eqnarray}\label{obs-1-9}
{{\partial u} \over {\partial t}} + {{\partial u} \over {\partial a}}+(k_{1}(x)u_{x})_{x}-\mu_{1}(t,a, x)u=\mu_{3}(t,a, x)v-\beta_{1}u(t,0,x)  && \text{ in } Q, \label{inter-1}\\
u(t,a, 1)=u(t,a, 0)=0 && \text{ on } (0,T)\times(0,A), \nonumber \\
\nonumber  u(T,a, x)=u_{T}(a, x) && \text{ in } Q_{A},\\
\nonumber  u(t,A, x)=0 && \text{ in } Q_{T},
\end{eqnarray}
and
\begin{eqnarray}\label{obs-1-10}
{{\partial v} \over {\partial t}} + {{\partial v} \over {\partial a}}+(k_{2}(x)v_{x})_{x}-\mu_{2}(t,a, x)v=-\beta_{2}v(t,0,x)  && \text{ in } Q, \label{inter-2}\\
v(t,a, 1)=v(t,a, 0)=0 && \text{ on } (0,T)\times(0,A), \nonumber \\
\nonumber  v(T,a, x)=v_{T}(a, x) && \text{ in } Q_{A},\\
\nonumber  v(t,A, x)=0 && \text{ in } Q_{T}.
\end{eqnarray}
Multiplying the first equations of \eqref{obser-1} and \eqref{obser-2} respectively by $\widetilde{u}$ and $\widetilde{v}$  and integrating
by parts on $Q_{t}=(0,t)\times(0,A)\times(0,1)$ one obtains
\begin{eqnarray}\label{obser-3}
\nonumber &&\frac{1}{2}\int_{Q_{A}}u^{2}(0,a,x)dadx
+\frac{1}{2}\int_{0}^{1}\int_{0}^{t}\widetilde{u}^{2}(\tau,0,x)d\tau dx
\\ \nonumber &&+\kappa\int_{0}^{1}\int_{0}^{A}\int_{0}^{t}\widetilde{u}^{2}(\tau,a,x)d\tau dadx \leq
\frac{\|\beta_{1}\|_{\infty}^{2}+1}{4\epsilon^{'}}\int_{0}^{1}\int_{0}^{A}\int_{0}^{t}\widetilde{u}^{2}(\tau,a,x)d\tau dadx
\\ &&+\epsilon^{'}A\int_{0}^{1}\int_{0}^{t}\widetilde{u}^{2}(\tau,0,x)d\tau dx+\epsilon^{'}\int_{Q_{t}}\mu_{3}^{2}\widetilde{v}^{2}d\tau dadx+\frac{1}{2}\int_{Q_{A}}\widetilde{u}^{2}(t,a,x)dadx.
\end{eqnarray}
and
\begin{eqnarray}\label{obser-4}
\nonumber &&\frac{1}{2}\int_{Q_{A}}v^{2}(0,a,x)dadx
+\frac{1}{2}\int_{0}^{1}\int_{0}^{t}\widetilde{v}^{2}(\tau,0,x)d\tau dx
\\ \nonumber &&+\kappa\int_{0}^{1}\int_{0}^{A}\int_{0}^{t}\widetilde{v}^{2}(\tau,a,x)d\tau dadx \leq
\frac{\|\beta_{2}\|_{\infty}^{2}+1}{4\epsilon^{'}}\int_{0}^{1}\int_{0}^{A}\int_{0}^{t}\widetilde{v}^{2}(\tau,a,x)d\tau dadx
\\ &&+\epsilon^{'}A\int_{0}^{1}\int_{0}^{t}\widetilde{v}^{2}(\tau,0,x)d\tau dx+\frac{1}{2}\int_{Q_{A}}\widetilde{v}^{2}(t,a,x)dadx.
\end{eqnarray}
Summing \eqref{obser-3} and \eqref{obser-4} side by side and taking  $\kappa=\max(\frac{\|\beta_{1}\|_{\infty}^{2}+1}{4\epsilon^{'}}, \frac{\|\beta_{2}\|_{\infty}^{2}+1}{4\epsilon^{'}}+\epsilon^{'}\|\mu_{3}\|_{\infty}^{2})$ and $\epsilon^{'}<\frac{1}{2A}$, on gets
\begin{eqnarray}\label{obs-1-5}
\int_{Q_{A}}u^{2}(0,a,x)dadx+\int_{Q_{A}}v^{2}(0,a,x)dadx\leq \int_{Q_{A}}\widetilde{u}^{2}(t,a,x)dadx +\int_{Q_{A}}\widetilde{v}^{2}(t,a,x)dadx.
\end{eqnarray}
Arguing as in \cite{Ain3}
 and integrating over $(\frac{T}{4}, \frac{3T}{4})$ we conclude
\begin{eqnarray}\label{obs-1-6}
\nonumber &&\int_{Q_{A}}u^{2}(0,a,x)dadx+\int_{Q_{A}}v^{2}(0,a,x)dadx\leq
C_{12}e^{2\kappa T}\left(\int_{0}^{1}\int_{0}^{\delta}u_{T}^{2}(a,x)dadx+\int_{0}^{1}\int_{0}^{\delta}v_{T}^{2}(a,x)dadx\right)
\\ && +\frac{2e^{2\kappa T}}{T}\left(\int_{0}^{1}\int_{\delta}^{A}\int_{\frac{T}{4}}^{\frac{3T}{4}}u^{2}(t,a,x)dtdadx
+\int_{0}^{1}\int_{\delta}^{A}\int_{\frac{T}{4}}^{\frac{3T}{4}}v^{2}(t,a,x)\right)dtdadx.
\end{eqnarray}
Hence, Hardy-Poincar\'{e} inequality and the definitions of $\varphi_{i}, i=1,2$ stated in \eqref{479} lead to
\begin{eqnarray}\label{obs-1-7}
\nonumber &&\int_{Q_{A}}u^{2}(0,a,x)dadx+\int_{Q_{A}}v^{2}(0,a,x)dadx\leq
C_{12}e^{2\kappa T}\left(\int_{0}^{1}\int_{0}^{\delta}u_{T}^{2}(a,x)dadx+\int_{0}^{1}\int_{0}^{\delta}v_{T}^{2}(a,x)dadx\right)
\\ \nonumber && +C_{\delta}^{13}\left(\int_{0}^{1}\int_{\delta}^{A}\int_{\frac{T}{4}}^{\frac{3T}{4}}s\Theta k_{1}(x) u^{2}(t,a,x)e^{2s\varphi_{1}}dtdadx+\int_{0}^{1}\int_{\delta}^{A}\int_{\frac{T}{4}}^{\frac{3T}{4}}s\Theta k_{2}(x) v^{2}(t,a,x)e^{2s\varphi_{2}}dtdadx\right).
\end{eqnarray}
 Finally, using the Carleman estimate \eqref{Car} we deduce the observability inequality \eqref{486}.
and then the proof is finished.
\endproof

Now, obtaining our observability inequality, following a standard argument, we are now ready to prove our
main result.
\begin{theorem}\label{487}
Assume that \eqref{473} and \eqref{474} are verified. Let $A>0$ and $T>0$ be given such that $T\in(0,\delta)$ with $\delta\in(0,A)$ small enough. Then, for all $(y_{0}, p_{0}) \in L^{2}(Q_{A})\times L^{2}(Q_{A})$,
there exists a control $\vartheta \in L^{2}(q)$ such that the associated solution of \eqref{470}  verifies
\begin{equation}\label{nul-contr-1}
\begin{cases}
&y(T,a, x)=0, \quad  \text{ a.e. in } (\delta, A)\times(0,1),\\
&p(T,a, x)=0, \quad   \text{ a.e in } (\delta, A)\times(0,1).
\end{cases}
\end{equation}
\end{theorem}

\proof

Let $ \varepsilon>0$ and consider the following cost function
 $$J_{\varepsilon}(\vartheta_{1},\vartheta_{2})=\frac{1}{2\varepsilon}\int_{0}^{1}\int_{\delta}^{A}(y^{2}(T,a, x)+p^{2}(T,a,x))dadx
+\frac{1}{2}\int_{q}\vartheta^{2}(t,a,x)dtdadx.$$
We can prove that $J_{\varepsilon}$ is continuous, convex and coercive. Then, it admits at least
one minimizer $\vartheta_{\varepsilon}$ and we have
\begin{eqnarray}\label{488}
\nonumber \vartheta_{\varepsilon}=-u_{\varepsilon}(t,a,x)\chi_{\omega}(x) \hspace{0.25cm} \text{ in } Q,\\
\end{eqnarray}

with $u_{\varepsilon}$ is the solution of the following system
\begin{eqnarray}\label{489}
{{\partial u_{\varepsilon}} \over {\partial t}} + {{\partial u_{\varepsilon}} \over {\partial a}}+(k_{1}(x)(u_{\varepsilon})_{x})_{x}-\mu_{1}(t,a, x)u_{\varepsilon}-\mu_{3}v_{\varepsilon} =-\beta_{1} u_{\varepsilon}(t,0,x) && \text{ in } Q,\\
\nonumber  u_{\varepsilon}(t,a, 1)=u_{\varepsilon}(t,a, 0)=0 &&  \text{ on }(0,T)\times (0,A),\\
\nonumber  u_{\varepsilon}(T,a, x)=\frac{1}{\varepsilon}y_{\varepsilon}(T,a, x)\chi_{(\delta, A)}(a) && \text{ in }Q_{A},\\
\nonumber  u_{\varepsilon}(t,A, x)=0 && \text{ in } Q_{T},
\end{eqnarray}
where $v_{\varepsilon}$ is the solution of
\begin{eqnarray}\label{490}
{{\partial v_{\varepsilon}} \over {\partial t}} + {{\partial v_{\varepsilon}} \over {\partial a}}+(k_{2}(x)(v_{\varepsilon})_{x})_{x}-\mu_{2}(t,a, x)v_{\varepsilon}=-\beta_{2} v_{\varepsilon}(t,0,x) && \text{ in } Q,\\
\nonumber  v_{\varepsilon}(t,a, 1)=v_{\varepsilon}(t,a, 0)=0   && \text{ on }(0,T)\times (0,A),\\
\nonumber  v_{\varepsilon}(T,a, x)=\frac{1}{\varepsilon}p_{\varepsilon}(T,a, x)\chi_{(\delta, A)}(a) && \text{ in }Q_{A},\\
\nonumber  v_{\varepsilon}(t,A, x)=0 && \text{ in } Q_{T},
\end{eqnarray}
 and $(y_{\varepsilon}, p_{\varepsilon})$ is the solution of the system \eqref{470} associated to the control
 $\vartheta_{\varepsilon}$. \\Multiplying the first equation of \eqref{489} by $y_{\varepsilon}$ and the second equation
 of \eqref{470} by $v_{\varepsilon}$, integrating over $Q$, using \eqref{488} and the Young inequality we obtain
\begin{eqnarray}\label{491}
\nonumber & &
\frac{1}{\varepsilon}\int_{0}^{1}\int_{\delta}^{A}(y_{\varepsilon}^{2}(T,a, x)+p_{\varepsilon}^{2}(T,a,x))dadx+\int_{q}\vartheta^{2}_{\varepsilon}(t,a,x) dtdadx\\ \nonumber&=&\int_{Q_{A}}(y_{0}(a,x)u_{\varepsilon}( 0, a, x)+p_{0}(a,x)v_{\varepsilon}(0, a, x))dadx \\ \nonumber& & \leq \frac{1}{4C_{\delta}}\int_{Q_{A}}(u_{\varepsilon}^{2}( 0, a, x)+v_{\varepsilon}^{2}( 0, a, x))dadx +C_{\delta}\int_{Q_{A}}(y_{0}^{2}(a,x)+p_{0}^{2}(a,x))dadx,
\end{eqnarray}
with $C_{\delta}$ is the constant of the observability inequality  \eqref{486}. Hence, using relation \eqref{488}, the observability inequality leads to
\begin{eqnarray}\label{491}
\nonumber & &\frac{1}{\varepsilon}\int_{0}^{1}\int_{\delta}^{A}(y_{\varepsilon}^{2}(T,a, x)+p_{\varepsilon}^{2}(T,a,x))dadx
+\frac{3}{4}\int_{q}\vartheta^{2}_{\varepsilon}(t,a,x) dtdadx \\ && \leq C_{\delta}\int_{Q_{A}}(y_{0}^{2}(a,x)+p_{0}^{2}(a,x))dadx.
\end{eqnarray}
Hence, it follows that
\begin{equation}\label{492}
\begin{cases}
\int_{0}^{1}\int_{\delta}^{A}y_{\varepsilon}^{2}(T,a, x)dadx \leq C_{\delta}\varepsilon\int_{Q_{A}}(y_{0}^{2}(a,x)+p_{0}(a,x))dadx,\\
\int_{0}^{1}\int_{\delta}^{A}p_{\varepsilon}^{2}(T,a, x)dadx \leq C_{\delta}\varepsilon\int_{Q_{A}}(y_{0}^{2}(a,x)+p_{0}(a,x))dadx,\\
\int_{q}\vartheta^{2}_{\varepsilon}(t,a,x) dtdadx \leq \frac{4C_{\delta}}{3}\int_{Q_{A}}(y_{0}^{2}(a,x)+p_{0}(a,x))dadx.\\
\end{cases}
\end{equation}
 Then, we can extract two subsequences of $(y_{\varepsilon},p_{\varepsilon})$ and $\vartheta_{\varepsilon}$
 denoted also by $\vartheta_{\varepsilon}$ and$(y_{\varepsilon},p_{\varepsilon})$ that converge weakly
 towards $\vartheta$ and $(y,p)$ in $L^{2}(q)$ and $L^{2}((0, T)\times(0, A); H^{1}_{k_{1}}(0, 1)\times H^{1}_{k_{2}}(0, 1))$
 respectively. Now, by a variational technic,  we prove that $(y,p)$ is a solution of  \eqref{470}  corresponding to the
controls $\vartheta$  and, by the first and second estimates of \eqref{492}, $(y,p)$  satisfies \eqref{471}.
\endproof
\section{Appendix}\label{last-section-5}
As is mentioned in the introduction, this section is devoted to the proofs of some intermediate results useful to show
the Carleman type inequality \eqref{Car}. Firstly, we begin by the Caccioppoli's inequality stated in the following lemma
\begin{lemma}\label{495} Let $\omega^{'}$ be a subset of $\omega$ such that $\omega^{'}\subset\subset\omega$. Then, there exists a positive constant $C$ such that
\begin{equation}\label{Caccio-1}
\int_{\omega^{'}}\int_{0}^{A}\int_{0}^{T} (u_{x}^{2}+v_{x}^{2})e^{2s\varphi_{i}} dtdadx \leq C \left(\int_{q}s^{2}\Theta^{2}(u^{2}+v^{2})e^{2s\varphi_{i}} dtdadx+\int_{q}(h_{1}^{2}+h_{2}^{2})e^{2s\varphi_{i}} dtdadx\right),
\end{equation}
 where $(u,v)$ is the solution of \eqref{478} and the weight functions $\varphi_{i}, i=1,2$ are defined by \eqref{479}.
\end{lemma}
\proof
The proof of this result is similar to the one of \cite[Lemma 5.1]{man}. Indeed, consider the cut-off function $\zeta$
defined by
\begin{equation}\label{496}
\left\{
\begin{array}{l}
0\leq\zeta(x)\leq1, \hspace{0.5cm} x \in \mathbb{R},\\
\zeta(x)=0, \hspace{0.5cm}   x<x_1 \text{ and } x>x_2,\\
\zeta(x)=1, \hspace{0.5cm}   x \in \omega^{'},\\
\end{array}
\right.
\end{equation}
For $(u,v)$ solution of \eqref{478} one has
\begin{eqnarray}\label{131}
\nonumber 0&=&\int_{0}^{T}\frac{d}{dt}\left[\int_{0}^{1}\int_{0}^{A}
\zeta^{2}e^{2s\varphi_{i}}(u^{2}+v^{2})dadx\right]dt
\\ \nonumber &=&2s \int_{0}^{1}\int_{0}^{A}
\int_{0}^{T}\zeta^{2}(\varphi_{i})_{t}(u^{2}+v^{2})e^{2s\varphi_{i}}dtdadx
+2\int_{0}^{1}\int_{0}^{A}
\int_{0}^{T}\zeta^{2}ww_{t}
e^{2s\varphi_{i}}dtdadx
\\ \nonumber &=&2s\int_{0}^{1}\int_{0}^{A}
\int_{0}^{T}\zeta^{2}(\varphi_{i})_{t}(u^{2}+v^{2})
e^{2s\varphi_{i}}dtdadx\\ \nonumber & &+2\int_{0}^{1}\int_{0}^{A}\int_{0}^{T}\zeta^{2}u(-(k_{1}u_{x})_{x}-u_{a}+h_{1}+\mu_{1}u+\mu_{3}v)e^{2s\varphi_{i}}dtdadx\\ \nonumber & &+
2\int_{0}^{1}\int_{0}^{A}\int_{0}^{T}\zeta^{2}v(-(k_{2}v_{x})_{x}-v_{a}+h_{2}+\mu_{2}v)e^{2s\varphi_{i}}dtdadx.
\end{eqnarray}
Then, integrating by parts we obtain
\begin{eqnarray}\label{82}
\nonumber 2\int_{Q}\zeta^{2}(k_{1}u_{x}^{2}+k_{2}v_{x}^{2})e^{2s\varphi_{i}}dtdadx&=&-2s\int_{Q}\zeta^{2}(u^{2}+v^{2})\psi_{i}(\Theta_{a}+
\Theta_{t})e^{2s\varphi_{i}}dtdadx\\ \nonumber &&-2\int_{Q}\zeta^{2}(uh_{1}+vh_{2})e^{2s\varphi_{i}}dtdadx-2\int_{Q}\zeta^{2}(\mu_{1}u^{2}+\mu_{2}v^{2})e^{2s\varphi_{i}}dtdadx
\\ \nonumber&&+\int_{Q}(k_{1}(\zeta^{2}e^{2s\varphi_{i}})_{x})_{x}u^{2}dtdadx+\int_{Q}(k_{2}(\zeta^{2}e^{2s\varphi_{i}})_{x})_{x}v^{2}dtdadx
\\ \nonumber &&-2\int_{Q}\zeta^{2}\mu_{3}uve^{2s\varphi_{i}}dtdadx.
\end{eqnarray}
On the other hand, by the definitions of $\zeta$, $\psi$ and $\Theta$, using Young inequality and taking $s$ quite large
there is a constant  $c$ such that
\begin{eqnarray}\label{88}
\nonumber&&2\int_{Q}\zeta^{2}(k_{1}u_{x}^{2}+k_{2}v_{x}^{2})e^{2s\varphi_{i}}dtdadx\geq 2\min(\min_{x\in\omega^{'}}k_{1}(x),\min_{x\in\omega^{'}}k_{2}(x))
\int_{\omega^{'}}\int_{0}^{A}\int_{0}^{T}(u_{x}^{2}+v_{x}^{2})e^{2s\varphi_{i}}dtdadx,\\
\nonumber&&\int_{Q}(k_{1}(\zeta^{2}e^{2s\varphi_{i}})_{x})_{x}u^{2}dtdadx\leq c\int_{\omega}\int_{0}^{A}\int_{0}^{T}s^{2}\Theta^{2}u^{2}e^{2s\varphi_{i}}dtdadx,\\
\nonumber&&\int_{Q}(k_{2}(\zeta^{2}e^{2s\varphi_{i}})_{x})_{x}v^{2}dtdadx\leq c\int_{\omega}\int_{0}^{A}\int_{0}^{T}s^{2}\Theta^{2}v^{2}e^{2s\varphi_{i}}dtdadx,\\
\nonumber&&-2s\int_{Q}\zeta^{2}(u^{2}+v^{2})\psi_{i}(\Theta_{a}+\Theta_{t})e^{2s\varphi_{i}}dtdadx\leq c\int_{\omega}\int_{0}^{A}\int_{0}^{T}s^{2}\Theta^{2}(u^{2}+v^{2})e^{2s\varphi_{i}}dtdadx,\\
\nonumber&&-2\int_{Q}\zeta^{2}uh_{1}e^{2s\varphi_{i}}dtdadx\leq c\left(\int_{\omega}\int_{0}^{A}\int_{0}^{T}s^{2}\Theta^{2}u^{2}e^{2s\varphi_{i}}dtdadx+\int_{\omega}\int_{0}^{A}\int_{0}^{T}h_{1}^{2}e^{2s\varphi_{i}}dtdadx\right),\\
\nonumber&&-2\int_{Q}\zeta^{2}vh_{2}e^{2s\varphi_{i}}dtdadx\leq c\left(\int_{\omega}\int_{0}^{A}\int_{0}^{T}s^{2}\Theta^{2}v^{2}e^{2s\varphi_{i}}dtdadx+\int_{\omega}\int_{0}^{A}\int_{0}^{T}h_{2}^{2}e^{2s\varphi_{i}}dtdadx\right),\\
\nonumber&&-2\int_{Q}\zeta^{2}(\mu_{1} u^{2}+\mu_{2}v^{2})e^{2s\varphi_{i}}dtdadx \leq c\int_{\omega}\int_{0}^{A}\int_{0}^{T}s^{2}\Theta^{2}(u^{2}+v^{2})e^{2s\varphi_{i}}dtdadx,\\
\nonumber&&-2\int_{Q}\zeta^{2}\mu_{3}uve^{2s\varphi_{i}}dtdadx\leq c\int_{\omega}\int_{0}^{A}\int_{0}^{T}s^{2}\Theta^{2}(u^{2}+v^{2})e^{2s\varphi_{i}}dtdadx.
\end{eqnarray}
Combining all these inequalities, we can see that there is $ C>0$ such that
\begin{eqnarray*}
& &\int_{\omega^{'}}\int_{0}^{A}\int_{0}^{T}(u_{x}^{2}+v_{x}^{2})e^{2s\varphi_{i}}dtdadx\leq C\left(\int_{q}s^{2}\Theta^{2}(u^{2}+v^{2})e^{2s\varphi_{i}} dtdadx+\int_{q}(h_{1}^{2}+h_{2}^{2})
e^{2s\varphi_{i}}dtdadx\right).
\end{eqnarray*}
Thus, the proof is achieved.
\endproof
\begin{remark}
In Lemma \ref{495}, the set $\omega^{'}$ is chosen so that $0$ which is exactly the point of degeneracy of the dispersion coefficients $k_{1}$ and $k_{2}$ does not belong to $\overline{\omega^{'}}$. More generally, if the degeneracy occurs at a point $x_{0}\in (0,1)$, one must take $x_{0}$ out of $\overline{\omega^{'}}$ in the case of interior degeneracy to establish a Caccioppoli's type inequality (see \cite{genni} for more details in this context).
\end{remark}
We close this section by the following result
\begin{lemma}\label{Intervalle}
Assume that the conditions \eqref{parameters} hold. Then, $I=[\frac{k_{2}(1)(2-\gamma)(e^{2\kappa\|\sigma\|_{\infty}}-1)}{d_{2}k_{2}(1)(2-\gamma)-1}, \frac{4(e^{2\kappa\|\sigma\|_{\infty}}-e^{\kappa\|\sigma\|_{\infty}})}{3d_{2}})$ is not empty.
\end{lemma}
\proof
Indeed, one has
\begin{eqnarray}\label{interval}
 \nonumber \frac{4(e^{2\kappa\|\sigma\|_{\infty}}-e^{\kappa\|\sigma\|_{\infty}})}{3d_{2}}
 -\frac{k_{2}(1)(2-\gamma)(e^{2\kappa\|\sigma\|_{\infty}}-1)}{d_{2}k_{2}(1)(2-\gamma)-1}
\\ \nonumber=\frac{4(e^{2\kappa\|\sigma\|_{\infty}}-e^{\kappa\|\sigma\|_{\infty}})(d_{2}k_{2}(1)(2-\gamma)-1)-3d_{2}k_{2}(1)(2-\gamma)(e^{2\kappa\|\sigma\|_{\infty}}-1)}{3d_{2}(d_{2}k_{2}(1)(2-\gamma)-1)}
\\ \nonumber=\frac{e^{2\kappa\|\sigma\|_{\infty}}(d_{2}k_{2}(1)(2-\gamma)-4)-4e^{\kappa\|\sigma\|_{\infty}}(d_{2}k_{2}(1)(2-\gamma)-1)}{3d_{2}(d_{2}k_{2}(1)(2-\gamma)-1)}
+\frac{k_{2}(1)(2-\gamma)}{d_{2}k_{2}(1)(2-\gamma)-1}
\\ \nonumber=\frac{e^{\kappa\|\sigma\|_{\infty}}[e^{\kappa\|\sigma\|_{\infty}}(d_{2}k_{2}(1)(2-\gamma)-4)-4(d_{2}k_{2}(1)(2-\gamma)-1)]}{3d_{2}(d_{2}k_{2}(1)(2-\gamma)-1)}
+\frac{k_{2}(1)(2-\gamma)}{d_{2}k_{2}(1)(2-\gamma)-1}.
 \end{eqnarray}
 Using the fact that $d_{2}\geq \frac {5}{k_{2}(1)(2-\gamma)}$, we can conclude that $\frac{4(d_{2}k_{2}(1)(2-\gamma)-1)}{d_{2}k_{2}(1)(2-\gamma)-4)}\leq 16$.\\
Since $\kappa\geq\frac{4\ln2}{\|\sigma\|_{\infty}}$, then we have $e^{\kappa\|\sigma\|_{\infty}}\geq 16$. Therefore, the
previous difference is positive and subsequently $I\neq\emptyset$.
\endproof


\begin{thebibliography}{99}
\bibitem{Ain4}B. Ainseba, Corrigendum to "Exact and approximate controllability of the age and space population dynamics structured model [J. Math. Anal. Appl. 275 (2002), 562-574]", J. Math. Anal. Appl. 393 (2012),  328.
\bibitem{Ain3}B. Ainseba, Exact and approximate controllability of the age and space population dynamics
structured model,  J. Math. Anal. Appl. 275 (2002), 562-574.
\bibitem{Ain5} B. Ainseba and S. Anita, Internal stabilizability for a reaction-diffusion problem modelling a predator-prey system, Nonlinear analysis, 61 (2005), 491-501.
\bibitem{Ain2} B. Ainseba and S. Anita, Internal exact controllability of the linear population dynamics with diffusion, Electronic Journal of Differential Equations, 2004(2004), 1-11.
\bibitem{Ain1}  B. Ainseba and S. Anita, Local exact controllability of the age-dependent population dynamics with diffusion, Abstr. Appl. Anal. 6 (2001), 357-368.
\bibitem{ech} B. Ainseba, Y. Echarroudi and L. Maniar, Null controllability of a population dynamics with degenerate diffusion, Journal of Differential and Integral Equations, Vol. 26, Number 11/12 (2013), pp.1397-1410.
\bibitem{Haj1}E. M. Ait ben hassi, F. Ammar Khodja, A. Hajjaj and L. Maniar, Null controllability of degenerate parabolic cascade systems, Portugaliae Mathematica, 68 (2011), 345-367.
\bibitem{Bouss} F. Alabau-Boussouira, P. Cannarsa and G. Fragnelli, Carleman estimates for degenerate parabolic operators with applications to null controllability, J. evol.equ 6 (2006), 161-204.
\bibitem{Marcheva} V. Barbu,  M. Iannelli and  M. Martcheva,  On the controllability of the Lotka-McKendrick model of population dynamics,  J. Math. Anal. Appl. 253 (2001),  142-165.
\bibitem{genni} G. Fragnelli and D. Mugnai, Carleman estimates and observability inequalities for
parabolic equations with interior degeneracy, Advances in Nonlinear Analysis 08/2013; 2(4):339–378. DOI: 10.1515/anona-2013-0015.
\bibitem{cmp}
\newblock M. Campiti, G. Metafune, and D. Pallara,
\newblock \emph{Degenerate self-adjoint evolution equations on the unit interval},
\newblock Semigroup Forum. 57(1998), pp. 1-36.

\bibitem{Can1}
\newblock P. Cannarsa, P. Martinez and J. Vancostenoble,
\newblock \emph{Null controllability of degenerate heat equations,}
\newblock Adv. Differential Equations. 10(2005), pp. 153-190.
\bibitem{man} Y. Echarroudi and L. Maniar, Null controllability of a model in population dynamics, Electronic Journal of Differential Equations, 2014 (2014), No. 240, 1-20.
\bibitem{Fursikov} A. V. Fursikov and O. Yu. Imanuvilov, Controllability of Evolution Equations, Lecture Notes Series, vol. 34, Seoul National University Research Institute of Mathematics Global Analysis Research
Center, Seoul, 1996.

\bibitem{Traore} O. Traore,  Null controllability of a nonlinear population dynamics problem, Int. J. Math. Sci.  (2006),  1-20.
\bibitem{web}   G. F. Webb, Population models structured by age, size, and spatial position. Structured population models in biology and epidemiology, 1–49, Lecture Notes in Math. 1936, Springer, Berlin, 2008.
\bibitem{Zh2} C. Zhao, M. Wang and P. Zhao, Optimal control of harvesting for age-dependent predator-prey system, Mathematical and Computer Modelling, 42 (2005), 573-584.
%
%
\end{thebibliography}
\end{document}